\documentclass[preprint,3p,number,sort&compress]{elsarticle}

\DeclareUnicodeCharacter{2212}{\ensuremath{-}}

\usepackage{color}
\usepackage{amssymb}
 \usepackage{amsthm}
\usepackage{amsmath}
\usepackage{algorithm, algorithmic}
\usepackage{cases}
\usepackage[hidelinks]{hyperref}
\usepackage{multirow}
\usepackage{booktabs}
\usepackage{subfigure}
\usepackage{caption}
\usepackage{colortbl}
\usepackage{array}

\newtheorem{definition}{Definition}[section]
\newtheorem{thm}[definition]{Theorem}

\newtheorem{emp}[definition]{Experiment}

\makeatletter
\@addtoreset{equation}{section}
\makeatother

\everymath{\displaystyle}

\journal{Journal of Computational and Applied Mathematics}

\begin{document}

\begin{frontmatter}



\title{{\bf A tensor Alternating Anderson--Richardson method for solving multilinear systems with $ \mathcal{M} $-tensors}}
\author[ngy]{Jing Niu\corref{cor1}}
\cortext[cor1]{Corresponding author.}
	\ead{j-niu@na.nuap.nagoya-u.ac.jp}
\author[dl]{Lei Du}
	\ead{dulei@dlut.edu.cn}
\author[ngy]{Tomohiro Sogabe}
	\ead{sogabe@na.nuap.nagoya-u.ac.jp}
\author[ngy]{Shao-Liang Zhang}
	\ead{zhang@na.nuap.nagoya-u.ac.jp}

\address[ngy]{Department of Applied Physics, Nagoya University, Furo-cho, Chikusa-ku, Nagoya 464-8603, Japan}
\address[dl]{School of Mathematical Sciences, Dalian University of Technology, Dalian, Liaoning, 116024, PR China}

\begin{abstract}
	
It is well-known that a multilinear system with a nonsingular $ \mathcal{M} $-tensor and a positive right-hand side has a unique positive solution.
Tensor splitting methods \textcolor{black}{generalizing} the classical iterative methods for linear systems have been proposed for finding the unique positive solution.
The Alternating Anderson--Richardson (AAR) method is an effective method to accelerate the classical iterative methods.
In this study, we apply the idea of AAR for finding the unique positive solution quickly.
We first present a tensor Richardson method based on tensor regular splittings, then apply Anderson acceleration to the tensor Richardson method and derive a tensor \textcolor{black}{Anderson--Richardson} method, finally, we periodically employ the tensor Anderson--Richardson method within the tensor Richardson method and propose a tensor AAR method.
Numerical experiments show that the proposed method is effective in accelerating tensor splitting methods.

\end{abstract}


\end{frontmatter}

\section{Introduction}
\label{sec:sec1} 
We consider iterative methods for solving a multilinear system
\begin{equation} 
	\label{Eq1.1}   
	\mathcal{A} \boldsymbol{x}^{m-1} =\boldsymbol{b} 
\end{equation}
with a nonsingular $ \mathcal{M} $-tensor $ \mathcal{A} $ and a vector $ \boldsymbol{b} \in \mathbb{R}_{++}^{n} $. $ \mathbb{R}_{++}^{n} $ denotes the set of all positive vectors in $ \mathbb{R}^{n}$.
The operation $\mathcal{A} \boldsymbol{x}^{m-1}$ is defined by
$$ \mathcal{A} \boldsymbol{x}^{m-1} = \mathcal{A} \bar{\times} _{2} \boldsymbol{x} \bar{\times} _{3} \boldsymbol{x} \bar{\times} _{4} \cdots \bar{\times} _{m} \boldsymbol{x} \text{. } $$
Generally, $ \mathcal{A} \bar{\times} _{k} \boldsymbol{x} $ represents the $k$-mode product of a tensor $ \mathcal{A} \in \mathbb{R}^{n \times n \times \dots \times n} $ with a vector $ \boldsymbol{x} \in \mathbb{R}^{n}  $.
We review the specific definitions of $ \mathcal{M} $-tensor and $ \mathcal{A} \bar{\times} _{k} \boldsymbol{x} $ in Section \ref{subsec:sec2.1}.

Multilinear systems often appear in numerical partial differential equations \cite{r3}, data mining \cite{r1}, and tensor complementarity problems \cite{r2}.
Owing to these broad applications, researching how to solve multilinear systems effectively has attracted significant attention recently.

Many tensor splitting methods have been proposed for solving Eq. \eqref{Eq1.1}.
Ding and Wei \cite{r3} extended the classical Jacobi and Gauss--Seidel methods for linear systems to solve Eq. \eqref{Eq1.1}.
These extended methods are equivalent to solving diagonal and triangular tensor equations.
They proved that a multilinear system with a nonsingular $ \mathcal{M} $-tensor $\mathcal{A} $ and a positive $ \boldsymbol{b} $ has a unique positive solution.
When the coefficient tensor is symmetric (called a symmetric system for short),
they also proposed a Newton method.
Li et al. \cite{r5} proposed several classical methods by minimizing the approximation of symmetric systems.
They also proposed a Newton--Gauss--Seidel method by using multistep Gauss--Seidel iterations.
Liu et al. \cite{r6} generalized some tensor splitting methods based on tensor regular splittings.
\textcolor{black}{They \cite{p1} also provided some spectral radius comparisons between two tensor splitting methods, and proposed a preconditioned tensor splitting method. There are various types of preconditioning techniques that can accelerate the convergence of tensor splitting methods, for more details see \cite{p1,p2,p3,p4,p5}.}  

Some other tensor-type methods for solving Eq. \eqref{Eq1.1} also exist.
For Eq. \eqref{Eq1.1} with $ \mathcal{M} $-tensors,
Han \cite{r8} proposed a homotopy method.
Xie et al. \cite{r9} based on the rank-1 approximation of a symmetric coefficient tensor $ \mathcal{A} $ and proposed a method.
He et al. \cite{r10} proposed a Newton-type method by rewriting Eq. \eqref{Eq1.1} as a nonlinear system involving P-functions.
Liu et al. \cite{r11} proposed a sufficient descent nonlinear conjugate gradient method with an inexact line search.
Li et al. \cite{r13} extended a Newton method for solving Eq. \eqref{Eq1.1} with a nonnegative right-hand side.
\textcolor{black}{In addition, some methods have been proposed for solving Eq. \eqref{Eq1.1} with other structured tensors. 
	Lv and Ma \cite{ss1} proposed a Levenberg--Marquardt method for solving Eq. \eqref{Eq1.1} with semi-symmetric coefficient tensors.
	Wang et al. \cite{ss2} proposed a preconditioned AOR iterative methods for solving Eq. \eqref{Eq1.1} with $ \mathcal{H} $-tensors.
	They also proposed two neural network models for solving Eq. \eqref{Eq1.1} with nonsingular tensors in \cite{r12}.
	Jiang and Li \cite{ss4} proposed a new preconditoned AOR-type method for Eq. \eqref{Eq1.1} with $ \mathcal{H} $-tensors.
	Wang et al. \cite{ss3} proposed two randomized Kaczmarz-like methods for solving Eq. \eqref{Eq1.1} with nonsingular tensors.}

We consider accelerating the convergence of the existing tensor splitting methods for solving Eq. \eqref{Eq1.1}. Anderson acceleration \cite{r14} is an efficient technique to accelerate the convergence of linear and nonlinear fixed-point iterations.
\textcolor{black}{Anderson acceleration} is a multisecant method \cite{scf1} and is essentially equivalent to the Generalized Minimal Residual (GMRES) method if all mixing parameters are equal to 1 \cite{r15}. 
The convergence has been studied \cite{aa1,aa2}.
Anderson acceleration has wide applications in many fields.
For example, Anderson acceleration is used to solve two classes of transport equations \cite{pcd2} or for improving the convergence of Picard iterations in variably saturated flow modeling \cite{pcd1}.
It is also an efficient procedure to accelerate the convergence of self-consistent field iterations in electronic structure calculations\cite{scf1,scf2,scf3}.
\textcolor{black}{However, Anderson acceleration needs to solve a least-squares problem.
	The computational cost of solving the least-squares problem at each iteration is high.
	To reduce the computational cost, an alternating Anderson acceleration technique has been proposed, which periodically employs Anderson acceleration within fixed-point iterations. The alternating Anderson acceleration technique is used to accelerate fixed-point iterations in \cite{r16, r17, r18}.}
Pratapa et al. \cite{r16} periodically employed Anderson acceleration within the Jacobi method and proposed an Alternating Anderson Jacobi (AAJ) method.
Inspired by the AAJ and Scheduled Relaxation Jacobi (SRJ) \cite{srj} methods, Kong et al. \cite{r18} proposed an Alternating Anderson--SRJ method.
Pratapa et al. \cite{r17} generalized the AAJ method with preconditioners and proposed an Alternating Anderson--Richardson (AAR) method.
Lupo Pasini \cite{aar3} proposed an augmented AAR method.

\textcolor{black}{The Gauss--Seidel method in the tensor splitting method \cite{r3} is based on a triangular splitting of a tensor. A forward substitution algorithm has been proposed for solving lower triangular tensor equations. However, the computational cost of the forward substitution algorithm is large, because it needs to compute a spectral radius of a matrix. To reduce the computational cost, the tensor splitting method \cite{r6} transforms Eq. \eqref{Eq1.1} into a matrix problem by using tensor regular splittings. Inspired by employing alternating Anderson acceleration in iterative methods for solving linear systems, it is expected that applying the alternating Anderson acceleration technique can improve the convergence of the tensor splitting method \cite{r6} for solving Eq. \eqref{Eq1.1}.} We first present a tensor Richardson (TR) method, then apply Anderson acceleration in the TR method and derive a tensor Anderson--Richardson (TAR) method, and finally, we \textcolor{black}{propose} a tensor AAR (TAAR) method.

The rest of this paper is organized as follows. 
We present relevant definitions and lemmas about tensors as well as review three tensor splitting methods \cite{r3,r5,r6} in Section \ref{sec:sec2}.
Then, we recall the AAR method for solving linear systems and propose the TAAR method for solving Eq. \eqref{Eq1.1} in Section \ref{sec:sec3}. 
Numerical experiments are performed to verify the effectiveness of the proposed method in Section \ref{sec:sec4}.

\section{Preliminaries}
\label{sec:sec2} 
We first present some elementary definitions and lemmas of tensors and then recall three tensor splitting methods \cite{r3,r5,r6}.

\subsection{Notions}
\label{subsec:sec2.1} 

\begin{definition}
	\emph{(\cite{r19,r5}).  Let $\mathbb{R}\left(\mathbb{C}\right)$ be the real (complex) field. An $ m $th\textcolor{black}{-}order $ n $-dimensional tensor is
		$$ \mathcal{A}=\left(a_{ i_{1} i_{2} \dots i_{m} }\right)\text{,} \quad   a_{ i_{1} i_{2} \dots  i_{m} } \in \mathbb{R}\text{,} \quad  1 \leq i_{1}\text{, } i_{2}\text{,} \dots \text{, }i_{m} \leq n \text{.}$$
		$ \mathbb{R}^{\left[m \text{,}n\right]} $ denotes the set containing all these tensors. \\
		\indent The diagonal part  $ \mathcal{D}$ of  $ \mathcal{A} $ contains entries $ a_{ ii \dots  i}\text{, } i=1\text{, } 2\text{,} \dots \text{, } n $, 
		and the lower triangular part  $ \mathcal{L}$ of $ \mathcal{A} $ contains entries $ a_{ i_{1} i_{2} \dots i_{m} } $ with $ i_{1} =1\text{, } 2\text{,} \dots \text{, } n $ and $ i_{2} \text{, } i_{3}\text{, }\dots\text{, } i_{m} \leqslant i_{1} $.
		$ \mathcal{A} $ is an identity tensor if all diagonal entries are equal to 1 and other entries are 0, denoted by $ \mathcal{I} $.
		$ \mathcal{A} $ is a diagonal face tensor \textcolor{black}{if its} entries satisfy $ a_{ i_{1} i_{2} i_{3} \dots i_{m} }=0\text{, } i_{1} \neq i_{2}$, denoted by $ \widetilde{\mathcal{D} } $. 
		$ \mathcal{A} $ is a lower half tensor \textcolor{black}{if its} entries satisfy $ a_{ i_{1} i_{2} i_{3} \dots i_{m} }=0\text{, } \forall i_{2}  >i_{1}   $, denoted by $ \widetilde{\mathcal{L} } $. 
	}
\end{definition}

\begin{definition}
	\emph{(\cite{r19}).
		Let $  \mathcal{A} \in \mathbb{R}^{\left[m\text{,}n\right]} $, $ \boldsymbol{x} \in \mathbb{R}^{n} $, the $k$-mode product $ \mathcal{A} \bar{\times} _{k} \boldsymbol{x}  $ is an $ \left(m-1\right) $th\textcolor{black}{-}order tensor, elementwise,
		$$ \left( \mathcal{A} \bar{\times} _{k} \boldsymbol{x} \right) _{ i_{1} \dots  i_{k-1}  i_{k+1} \dots  i_{m}} = \sum \limits_{ i_{k}=1 }^{n} a_{ i_{1} \dots i_{k} \dots i_{m}} \boldsymbol{x}_{i_{k}}\text{.}$$
		$ \mathcal{A} \boldsymbol{x}^{m-1} $ is an $n$-dimensional vector defined by the $k$-mode product,
		\begin{align*}
			\left(\mathcal{A} \boldsymbol{x}^{m-1}\right)_{i} &= \left(\mathcal{A} \bar{\times} _{2} \boldsymbol{x} \bar{\times} _{3} \dots \bar{\times} _{m} \boldsymbol{x}\right)_{i} \\
			&= \sum \limits_{i_{2}\text{,}\dots\text{, }i_{m}=1}^{n} a_{ i i_{2} \dots i_{m}}  \boldsymbol{x}_{i_{2}} \dots \boldsymbol{x}_{i_{m}}\text{,} \quad  i=1\text{, } 2\text{,} \dots\text{, }n\text{. }
		\end{align*}
		$ \mathcal{A} \boldsymbol{x}^{m-2} $ is an $ n \times n $ matrix, elementwise,
		\begin{equation} \nonumber
			\left(\mathcal{A} \boldsymbol{x}^{m-2} \right)_{ij} = \left( \mathcal{A} \bar{\times} _{3} \boldsymbol{x} \bar{\times} _{4} \dots \bar{\times} _{m} \boldsymbol{x} \right)_{ij}\text{.}
		\end{equation}
	}
\end{definition}

\begin{definition}
	\label{def7}
	\emph{(\cite{r24}). 
		Let $ \mathcal{A} \in \mathbb{R}^{\left[m\text{,}n\right]}$, $ \mathcal{B} = \left(\mathcal{B}_{ i_{1} i_{2} \dots i_{k} }\right) \in \mathbb{R}^{\left[k\text{,}n\right]}$, 
		the product $ \mathcal{A} \mathcal{B}  $ is an $\left(m-1\right)\left(k-1\right) + 1$th\textcolor{black}{-}order tensor with entries
		$$ \left(\mathcal{A} \mathcal{B} \right)_{j \alpha_{2}\dots \alpha_{m}} =  \sum_{j_{2}\text{,}\dots\text{, }j_{m} =1}^{n} \left(a_{j j_{2}\dots j_{m}}\prod_{i=2}^{m}\mathcal{B}_{j_{i} \alpha_{i}}\right) \text{,}$$
		where $ j = 1\text{, } 2\text{,} \dots\text{, } n$, $ \alpha_{\textcolor{black}{i}} = \alpha_{\textcolor{black}{i}} ^{1} \alpha_{\textcolor{black}{i}} ^{2} \dots \alpha_{\textcolor{black}{i}} ^{k-1}$ with $ \alpha_{\textcolor{black}{i}} ^{1}\text{, } \alpha_{\textcolor{black}{i}} ^{2}\text{,} \dots\text{, } \alpha_{\textcolor{black}{i}} ^{k-1} \in \left \langle n \right \rangle   \text{, } i=2\text{, } 3\text{,} \dots\text{, } m\text{, } \left \langle n \right \rangle :=\{1\text{, }2\text{, } 3\text{,} \dots n\}$.}
\end{definition}

We consider a special case of Definition \ref{def7}.
Let $ A \in \mathbb{R}^{n \times n} $, the matrix-tensor product $ \mathcal{C} = A \mathcal{B} $ is a $ k $th\textcolor{black}{-}order tensor with entries
$$ \left(A \mathcal{B} \right)_{i_{1} i_{2}\dots i_{k}} =  \sum_{j_{2}=1 }^{n} \left(a_{i_{1} j_{2}}\mathcal{B}_{j_{2}i_{2}\dots i_{k}}\right) \text{,} \quad  i_{1}\text{, } i_{2}\text{,}\dots\text{, } i_{k} \in  \left \langle n \right \rangle  \text{.}$$

\begin{definition}
	\emph{(\cite{r20,r21}). 
		If $ \lambda \in \mathbb{C} $ and $ \boldsymbol{x} \in \mathbb{C}^{n}  \backslash \{ 0 \} $ satisfy
		$$ \mathcal{A} \boldsymbol{x}^{m-1} = \lambda \boldsymbol{x}^{[m-1]}\text{,} $$
		where $ \boldsymbol{x}^{\left[m-1\right]} =\left[\boldsymbol{x}_{1}^{m-1} \text{, } \boldsymbol{x}_{2}^{m-1}\text{,}\dots \text{, }\boldsymbol{x}_{n}^{m-1}\right]^{\top} $ and
		$\boldsymbol{x}_i^{m-1}$ represents the $(m-1)$\textcolor{black}{th} power of $\boldsymbol{x}_i$,
		then $ \lambda $ is an eigenvalue of $ \mathcal{A} $ and $ \boldsymbol{x} $ is the corresponding eigenvector. 
		The spectral radius $ \rho\left(\mathcal{A} \right) $ is defined by
		$$ \rho\left(\mathcal{A} \right)= \max\{|\lambda|:\lambda \text{ is an eigenvalue of } \mathcal{A}\}\text{.}  $$}
\end{definition}

\begin{definition}
	\emph{(\cite{r22}). 
		If a tensor $ \mathcal{A} \in \mathbb{R}^{\left[m\text{,}n\right]} $ satisfies 
		$$ \mathcal{A}=s \mathcal{I}-\mathcal{B} $$
		with  a nonnegative tensor  $ \mathcal{B} \in \mathbb{R}^{\left[m\text{,}n\right]} $ and a real number $ s >0 $, then $ \mathcal{A} $ is \textcolor{black}{a} $ \mathcal{Z} $-tensor. In particular, if $ s\geqslant \rho\left(\mathcal{B}\right) $, then $ \mathcal{A} $ is an $ \mathcal{M} $-tensor. If $ s> \rho\left(\mathcal{B}\right) $, $ \mathcal{A} $ is a nonsingular $ \mathcal{M} $-tensor.
		By default, $ \mathcal{A} $  below is a nonsingular $ \mathcal{M} $-tensor if not specified.}
\end{definition}

\begin{thm}
	\emph{(\cite{r22}, \textcolor{black}{Theorem 2 and Theorem 3}).
		\label{lma1}
		If $  \mathcal{A} \in \mathbb{R}^{[m\text{,}n]} $ is a $\mathcal{Z}$-tensor, the following conditons are equivalent
		\begin{enumerate}
			\item $ \mathcal{A} $ is a nonsingular $ \mathcal{M} $-tensor,
			\item There exists $ x>0 $ with $ \mathcal{A} x^{m-1} >0$.
		\end{enumerate}
	}
\end{thm}

Next, we present relevant definitions and lemmas about regular splittings.
\begin{definition}
	\emph{(\cite{r25}).
		Let $ \mathcal{A} \in \mathbb{R}^{\left[m\text{,}n\right]} $, $ \mathcal{B} \in \mathbb{R}^{\left[k\text{,}n\right]} $, if $ \mathcal{A} \mathcal{B}=\mathcal{I} $, then $ \mathcal{A} $ is an \textcolor{black}{$m$th-order} left inverse of $ \mathcal{B} $, and $ \mathcal{B} $ is \textcolor{black}{a $ k $th-order} right inverse of $ \mathcal{A} $.
		$ \mathcal{A} $ is a left (right)-nonsingular tensor or left (right)-invertible tensor if $ \mathcal{A} $ has \textcolor{black}{a $ k $th-order} left (right) inverse with $ k \geqslant 2 $.
	}
\end{definition}
\begin{definition}
	\emph{(\cite{r26}).
		Let $ R_{i}\left(\mathcal{A}\right) $ satisfy
		$$ R_{i}\left(\mathcal{A} \right) =\left( a_{ii_{2} \dots i_{m}} \right) _{i_{2} \dots i_{m}} ^{n} \in \mathbb{R}^{\left[m-1\text{,}n\right]}\text{,} \quad  1 \leqslant i \leqslant n\text{.} $$ 
		If all $ R_{i}\left(\mathcal{A}\right) $ are diagonal tensors, $ \mathcal{A} $ is row diagonal.
	}
\end{definition}

\begin{definition}
	\emph{(\cite{mm1}).
		The majorization matrix $ \mathit{M}\left(\mathcal{A}\right) $ of $ \mathcal{A} $ is an $ n \times n $ matrix with entries
		$$ \left(\mathit{M}\left(\mathcal{A}\right)\right) _{ij} := a_{ij \dots j}\text{,} \quad  i\text{, }j=1\text{, } 2\text{,} \dots \text{, } n. $$
		We use an example to easily understand $ \mathit{M}\left(\mathcal{A}\right) $. Let $ \mathcal{A} \in \mathbb{R}^{[3,3]} $ be \\
	}
\end{definition}
\begin{center} $ \mathcal{A}_{::1} =
	\left[
	\begin{matrix}
		1 & 2 &3 \\
		4 & 5 & 6 \\
		7 & 8 & 9
	\end{matrix}
	\right]
	$, \
	$ \mathcal{A}_{::2} =
	\left[
	\begin{matrix}
		10 & 11 &12 \\
		13& 14 & 15 \\
		16 & 17& 18
	\end{matrix}
	\right]
	$, \
	$ \mathcal{A}_{::3} =
	\left[
	\begin{matrix}
		19 & 20 &21 \\
		22& 23 & 24 \\
		25 & 26& 27
	\end{matrix}
	\right]
	$,
\end{center}
then $ \mathit{M}\left(\mathcal{A}\right) $ is 
$$ \mathit{M}\left(\mathcal{A}\right) =  \left[
\begin{matrix}
	1 & 11 &21 \\
	4& 14 & 24 \\
	7 & 17& 27
\end{matrix}
\right] \text{. } $$

\begin{thm}
	\emph{\textcolor{black}{(\cite{r6}, Lemma 3.6).
			\label{thm0}
			If $ \mathcal{A} $ is a nonsingular $ \mathcal{M} $-tensor, then $ \mathit{M}\left(\mathcal{A}\right) $  is a nonsingular $ \mathcal{M} $-matrix.}}
\end{thm}

\begin{thm}
	\emph{(\cite{r26}, \textcolor{black}{Proposition 5.1}).
		Let $ \mathcal{A} \in \mathbb{R}^{[m\text{,}n]} $, $ \mathcal{A} $ is row diagonal if and only if
		$ \mathcal{A} $ satisfies $ \mathcal{A}= \mathit{M}\left(\mathcal{A} \right) \mathcal{I} \text{. } $
	}
\end{thm}

\begin{thm}
	\emph{(\cite{r28}, \textcolor{black}{Theorem 3.1 and Corollary 3.3}).
		Let $ \mathcal{A} \in \mathbb{R}^{[m\text{,}n]} $, $ \mathcal{A} $ has \textcolor{black}{a} unique \textcolor{black}{2nd-order} left inverse $ \mathit{M}\left(\mathcal{A}\right)^{-1} $ 
		if and only if $ \mathcal{A} $ is row diagonal and $ \mathit{M}\left(\mathcal{A}\right) $ is nonsingular.
	}
\end{thm}

\begin{definition}
	\emph{(\cite{r6}).
		Let $ \mathcal{A}\text{, } \mathcal{E}\text{, and } \mathcal{F} \in \mathbb{R}^{\left[m\text{,}n\right]} $. $ \mathcal{A}=\mathcal{E}-\mathcal{F} $ is a regular splitting of $ \mathcal{A} $ 
		if $  \mathcal{F} \geqslant 0 $ and $ \mathcal{E}= \mathit{M}\left(\mathcal{E} \right) \mathcal{I} $ is left-nonsingular with $   \mathit{M}\left(\mathcal{E}\right)^{-1} \geqslant 0 $.
	}
\end{definition}


\begin{thm}
	\emph{(\cite{r3}, \textcolor{black}{Theorem 3.2}).
		\label{thm1}
		The multilinear system \eqref{Eq1.1} has a unique positive solution if $\mathcal{A}$ is a nonsingular $\mathcal{M}$-tensor and $\boldsymbol{b}$ is a positive vector.}
\end{thm}

\subsection{Tensor splitting methods}
\label{subsec:subsec2.2} 
Three tensor splitting methods \cite{r3,r5,r6} have been proposed for solving Eq. \eqref{Eq1.1}.
We denote these three methods into tensor splitting methods 1, 2, and 3 based on their inventors. 
Because the successive over relaxation (SOR)-type iterations of the tensor splitting methods need to choose proper relaxation parameters,
we do not compare the SOR-type iterations.
We review these three tensor splitting methods as follows.

\subsubsection{Tensor splitting method 1}
\label{subsubsec:sec2.1} 
Let $ \mathbb{R}_{++}^{n} $ denote the set of all positive real numbers, the positive solution set of Eq. \eqref{Eq1.1} is
$$ \left( \mathcal{A}^{-1}\boldsymbol{b} \right)_{++}:=\{ \boldsymbol{x}\in \mathbb{R}_{++}^{n}:  \mathcal{A} \boldsymbol{x}^{m-1}=\boldsymbol{b}\}\text{,} $$
\textcolor{black}{which has a unique element by Theorem \ref{thm1}. We denote the unique positive solution of Eq. \eqref{Eq1.1} as $  \mathcal{A}^{-1}_{++}\boldsymbol{b} $.}
According to the splittings of $\mathcal{A}  $, there are three iterative methods in \cite{r3}.
\begin{enumerate}[(1)]
	\item Jacobi (J1) method \\
	Let $ \mathcal{A}=\mathcal{D}-\mathcal{F} $, $\mathcal{D}$ is the diagonal part of $ \mathcal{A} $, the J1 method is
	$$\textcolor{black}{ \boldsymbol{x}_{k+1}= \mathcal{D} ^{-1}_{++} \left(\mathcal{F} \boldsymbol{x}_{k}^{m-1} + \boldsymbol{b}\right)  \text{,} \quad  k=0\text{, } 1\text{,}\dots\text{.} } $$
	\textcolor{black}{The solution $\boldsymbol{x}_{k+1}$ is in $ \left( \mathcal{D} ^{-1} \left(\mathcal{F} \boldsymbol{x}_{k}^{m-1} + \boldsymbol{b}\right) \right)_{++}$.} Since 
	$$ \mathcal{D} \boldsymbol{x}_{k+1}^{m-1} = \left[a_{11 \dots 1}\left(\boldsymbol{x}_{k+1}\right)_{1}^{m-1}\text{, } a_{22 \dots 2}\left(\boldsymbol{x}_{k+1}\right)_{2}^{m-1}\text{,}  \dots\text{, } a_{nn \dots n}\left(\boldsymbol{x}_{k+1}\right)_{n}^{m-1}\right]^{\top}\text{,}$$
	$ \boldsymbol{x}_{k+1} $ is
	$$ \left(\boldsymbol{x}_{k+1}\right)_{i} = \left(\frac{\left(\mathcal{F} \boldsymbol{x}_{k}^{m-1} + \boldsymbol{b}\right)_{i}}{a_{ii \dots i}}\right)^{\frac{1}{m-1}} \text{,} \quad  i= 1\text{, } 2\text{,} \dots\text{, } n\text{.} $$
	\item Gauss--Seidel (GS1) method \\
	Let $ \mathcal{A}=\mathcal{L}-\mathcal{F} $, $ \mathcal{L}$ is the lower triangular part of $ \mathcal{A} $, the GS1 method is
	$$\textcolor{black}{ \boldsymbol{x}_{k+1}= \mathcal{L} ^{-1}_{++} \left(\mathcal{F} \boldsymbol{x}_{k}^{m-1} + \boldsymbol{b}\right) \text{,} \quad  k=0\text{, } 1\text{,} \dots\text{.}}$$
	Ding and Wei \cite{r3} proposed a forward substitution algorithm for solving lower triangular tensor equations. 
	\item SOR-like method \\
	Ding and Wei also proposed a SOR-like method to accelerate the above methods. Let $\mathcal{A} = \mathcal{M}-\mathcal{N}$, the SOR-like method is
	$$ \textcolor{black}{ \boldsymbol{x}_{k+1}= (\mathcal{M}-\omega \mathcal{I}) ^{-1}_{++} \left((\mathcal{N}-\omega \mathcal{I}) \boldsymbol{x}_{k}^{m-1} + \boldsymbol{b}\right) \text{,} \quad  k=0\text{, } 1\text{,} \dots\text{,}} $$
	where $ \mathcal{M} $ is chosen as the diagonal part $ \mathcal{D} $ or the lower triangular part $ \mathcal{L} $ of $ \mathcal{A} $. The corresponding methods are denoted as J1\_SORlike and GS1\_SORlike. In \cite{r3}, the acceleration parameter $\omega $ is chosen as 
	$$ \omega = 0.35 \cdot \min _{i = 1 \text{, }2 \text{, } \dots \text{, } n} a_{ii \dots i} \text{.}$$
	\item Newton method for Eq. \eqref{Eq1.1} with Symmetric $\mathcal{M}$-tensors \\
	Solving Eq. \eqref{Eq1.1} with Symmetric $\mathcal{M}$-tensors is equivalent to solving the problem
	$$ \min_{\boldsymbol{x}\in \Omega}  \varphi(\boldsymbol{x}) := \frac{1}{m}\mathcal{A}\boldsymbol{x}^{m}-\boldsymbol{x}^{\top}\boldsymbol{b}\text{, }$$
	where $ \Omega=\{ x>0: \mathcal{A}\boldsymbol{x}^{m-1}>0 \} $.
	Ding and Wei \cite{r3} employed the Newton method. The iteration is $$ \boldsymbol{x}_{k+1} = M_{k}^{-1} (\frac{m-2}{m-1}\mathcal{A}\boldsymbol{x}_{k}^{m-1}+\frac{1}{m-1}\boldsymbol{b}) \text{,} $$
	where $  M_{k} = \mathcal{A} \boldsymbol{x}_{k}^{m-2} $ is a matrix.
\end{enumerate}

\subsubsection{Tensor splitting method 2}
\label{subsec:sec2.2} 
Because general tensors can be partially symmetrized, Li et al. \cite{r5} considered symmetric multilinear systems and proposed two iterative methods based on splittings of $\mathcal{A}$.
\begin{enumerate}[(1)]
	\item Jacobi (J2) method \\ [0.1em]
	Let $ \mathcal{A}=\widetilde{\mathcal{D}} -\mathcal{F} $ with the diagonal face part $ \widetilde{\mathcal{D} } $ of $ \mathcal{A} $, the J2 method is
	$$ \boldsymbol{x}_{k+1} = \boldsymbol{x}_{k} + \frac{1}{m-1} \left(\widetilde{\mathcal{D}}  \boldsymbol{x}_{k}^{m-2}\right)^{-1} \left( \boldsymbol{b}- \mathcal{A} \boldsymbol{x}_{k}^{m-1} \right) \text{,} \quad  k=0\text{, } 1\text{,} \dots\text{,}$$
	where $ \widetilde{\mathcal{D}}  \boldsymbol{x}_{k}^{m-2} $ is a diagonal matrix.
	\item Gauss--Seidel (GS2) method \\ [0.1em]
	Let $ \mathcal{A}=\widetilde{\mathcal{L}} -\mathcal{F} $ with the lower half part  $ \widetilde{\mathcal{L} } $ of $ \mathcal{A} $, the GS2 method is
	$$ \boldsymbol{x}_{k+1} = \boldsymbol{x}_{k} + \frac{1}{m-1} \left(\widetilde{\mathcal{L}}  \boldsymbol{x}_{k}^{m-2}\right)^{-1} \left( \boldsymbol{b}- \mathcal{A} \boldsymbol{x}_{k}^{m-1} \right) \text{,} \quad  k=0\text{, } 1\text{,} \dots\text{,}$$
	where $ \widetilde{\mathcal{L}}  \boldsymbol{x}_{k}^{m-2} $ is a lower triangular matrix.
\end{enumerate}

\subsubsection{Tensor splitting method 3}
\label{subsec:sec2.3} 
Liu et al. \cite{r6} proposed three iterative methods based on regular splittings of $\mathcal{A}$:
\begin{enumerate}[(1)]
	\item Jacobi (J3) method \\
	Consider a regular splitting of $ \mathcal{A} $: 
	$ \mathcal{A}= D\left( \mathit{M}\left(\mathcal{A}\right)\right) \mathcal{I} - \mathcal{F}\text{,}$ where $ D\left( \mathit{M}\left(\mathcal{A}\right)\right) $ is the diagonal part of $ \mathit{M}\left(\mathcal{A}\right) $. The J3 method is
	$$ \boldsymbol{x}_{k+1} = \left[ D\left( \mathit{M}\left(\mathcal{A}\right)\right)^{-1} \left(\boldsymbol{b}+ \mathcal{F} \boldsymbol{x}_{k}^{m-1}\right) \right] ^{[\frac{1}{m-1}]} \text{,} \quad  k=0\text{, } 1\text{,} \dots\text{,}$$
	where $  \mathit{M}\left(\mathcal{A}\right)_{ij} = a _{ij \dots j}\text{, }  i\text{, }j=1\text{, } 2\text{,} \dots \text{, } n $ and $ \boldsymbol{x}^{[\frac{1}{m-1}]} =\left[\boldsymbol{x}^{\frac{1}{m-1}}_{1}\text{, } \boldsymbol{x}^{\frac{1}{m-1}}_{2}\text{,} \dots\text{, }  \boldsymbol{x}^{\frac{1}{m-1}}_{n} \right]^{ \textcolor{black}{\top}} \in \mathbb{R}^{n}$.
	\item Gauss--Seidel (GS3) method \\ 
	Consider a regular splitting of $ \mathcal{A} $: 
	$ \mathcal{A}= L\left( \mathit{M}\left(\mathcal{A}\right)\right) \mathcal{I} - \mathcal{F} $, where $ L \left( \mathit{M}\left(\mathcal{A}\right)\right) $ is the lower triangular part of $ \mathit{M}\left(\mathcal{A}\right) $, the GS3 method is
	$$ \boldsymbol{x}_{k+1} = \left[ L\left( \mathit{M}\left(\mathcal{A}\right)\right)^{-1} \left(\boldsymbol{b}+ \mathcal{F} \boldsymbol{x}_{k}^{m-1}\right) \right] ^{[\frac{1}{m-1}]}\text{,} \quad  k=0\text{, } 1\text{,} \dots \text{,}$$
	where $ L\left( \mathit{M}\left(\mathcal{A}\right)\right)^{-1} \geqslant 0$ and $\mathcal{F}  \geqslant 0$. 
	\item FULLM method \\
	Let $ \mathcal{A}= \mathit{M}\left(\mathcal{A}\right) \mathcal{I} - \mathcal{F} $. 
	\textcolor{black}{The matrix $ \mathit{M}\left(\mathcal{A}\right)  $ is a nonsingular $ \mathcal{M} $-matrix by Theorem \ref{thm0}. 
		The matrix $ \mathit{M}\left(\mathcal{A}\right)  $ is inverse-positive \cite{np} and $ \mathit{M}\left(\mathcal{A}\right) ^{-1} \ge 0$.
		The tensor $ \mathcal{F} =  \mathit{M}\left(\mathcal{A}\right) \mathcal{I}-\mathcal{A} \ge 0$. So $ \mathcal{A}= \mathit{M}\left(\mathcal{A}\right) \mathcal{I} - \mathcal{F} $ is a regular splitting of $ \mathcal{A}$.
		Liu et al. \cite{r6} proposed a FULLM method based on the regular splitting.} The FULLM method is
	$$ \boldsymbol{x}_{k+1} = \left[ \mathit{M}\left(\mathcal{A}\right)^{-1} \left(\boldsymbol{b}+ \mathcal{F} \boldsymbol{x}_{k}^{m-1}\right) \right]^{[\frac{1}{m-1}]}\text{,} \quad  k=0\text{, } 1\text{,} \dots \text{,}$$
	where $\mathcal{F}  \geqslant 0$ and $\mathit{M}\left(\mathcal{A}\right) \mathcal{I} $ is a left-nonsingular tensor with $\mathit{M}\left(\mathcal{A}\right)^{-1} \geqslant 0$. 
\end{enumerate}

\section{Alternating Anderson--Richardson method} 
\label{sec:sec3} 
We first review the AAR method for solving linear systems,
then \textcolor{black}{propose} a TAAR method for solving Eq. \eqref{Eq1.1}.

\subsection{Alternating Anderson--Richardson method} 
\label{subsec:subsec3.1}
\textcolor{black}{Consider} a linear system
\begin{equation} \label{Eq3.1}   A\boldsymbol{x}=\boldsymbol{b}\text{,}  \end{equation}
where $ A \in \mathbb{R}^{n \times n} $ is a nonsingular matrix and $ \boldsymbol{b} \in \mathbb{R}^{n} \backslash \{0\} $. 

\subsubsection{Richardson method}
Let $ A= D-N $, where $ D $ is the diagonal part of $ A $. 
The weighted Jacobi method is 
$$ \boldsymbol{x}_{k+1} = \boldsymbol{x}_{k} + \omega_{k} D^{-1} \left(\boldsymbol{b}-A\boldsymbol{x}_{k}\right)\text{,}  $$
where $ D $ is a preconditioner. 
Using a different preconditioner $ M $, the Richardson method is 
\begin{equation} \label{Eq3.2}   \boldsymbol{x}_{k+1} = \boldsymbol{x}_{k} + \omega_{k} \boldsymbol{r}_{k} \text{,} \end{equation}
where the residual is $ \boldsymbol{r}_{k}=  M^{-1} \left( \boldsymbol{b}-A \boldsymbol{x}_{k} \right) $. 
The relaxation parameter $ \omega_{k}  $ is chosen by minimizing the quasi\textcolor{black}{-}residual as follows:
\begin{align*}
	\omega_{k} &= \arg\min  \Vert  \boldsymbol{b}-A \boldsymbol{x}_{k+1}  \Vert _{2} \\
	&= \arg\min  \Vert \boldsymbol{b}-A \boldsymbol{x}_{k} -\omega_{k} A \boldsymbol{r}_{k} \Vert _{2}\text{.}
\end{align*}
The optimal solution $ \omega_{k} $ satisfies
$$ \omega_{k} = \frac{\left(\boldsymbol{b}-A\boldsymbol{x}_{k}\text{, } A \boldsymbol{r}_{k}\right)}{\left(A\boldsymbol{r}_{k}\text{, } A\boldsymbol{r}_{k}\right)} \text{.}$$
Two types of preconditioner $ M $ are given in \cite{r17}.
\begin{enumerate}[(1)]
	\item Jacobi  preconditioner \\
	The AAJ method is recovered by using preconditioner $ D \left(A\right) $, which is the diagonal part of $ A $.
	\item ILU$(0)$ preconditioner\\
	For a large sparse matrix $ A $, the LU decomposition $ A = LU $ may cause the sparsity of $ L\text{, } U $ to be less than matrix $A$.
	This phenomenon is called fill-in, which increases the storage \textcolor{black}{cost}.
	To fix the fill-in problem, we can eliminate some nonzero elements of $ L\text{, } U $ and produce $ \widetilde{L}\text{, } \widetilde{U} $ that the sparsity pattern is the same as $ A $.
	The product $ \widetilde{L} \widetilde{U} $ is the incomplete LU factorization with no fill-in, denoted by ILU$(0)$.
\end{enumerate}
Let $ \boldsymbol{x}^{*} =A^{-1}\boldsymbol{b} $, the error is
\begin{align*}
	\boldsymbol{e}_{k+1} &= \boldsymbol{x}_{k+1}- \boldsymbol{x}^{*} \\  
	&= \left(I-\omega_{k} M^{-1}A\right)\boldsymbol{e}_{k}\text{.} 
\end{align*} 
The residual is $ \boldsymbol{r}_{k+1} = \left(I-\omega_{k} M^{-1}A\right)\boldsymbol{r}_{k} $
and the convergence requires $ \rho \left( I-\omega_{k} M^{-1}A \right) <1 $.

\subsubsection{Anderson--Richardson method}
\label{subsubsec:subsubsec3.1.2}
The Anderson--Richardson (AR) method employs Anderson acceleration in the Richardson method.
The AR method involves two steps:

The first step uses parameters $ \gamma_{1}\text{, } \gamma_{2}\text{,} \dots\text{, } \gamma_{q} \in \mathbb{R} $ to correct $ \boldsymbol{x}_{k} $ in Eq. \eqref{Eq3.2} as follows:
\begin{equation}
	\label{Eq3.3}
	\overline{\boldsymbol{x}}_{k} =\boldsymbol{x}_{k}- \sum_{j=1}^{q} \gamma_{j} \left( \boldsymbol{x}_{k-q+j} - \boldsymbol{x}_{k-q+j-1}\right)\text{.}
\end{equation} 
Let $ \boldsymbol{\Gamma}_{q}  \in \mathbb{R}^{q}  $, define $ X_{k}\text{ and } R_{k}  \in \mathbb{R}^{n \times q} $ as the iteration and residual histories of the $k$th iteration:
\begin{align*}
	\boldsymbol{\Gamma}_{q} &= \left[\gamma_{1}\text{, } \gamma_{2}\text{,} \dots\text{, } \gamma_{q} \right]^{\top}\text{,}\\
	X_{k} &= \left[\left(\boldsymbol{x}_{k-q+1}-\boldsymbol{x}_{k-q}\right)\text{, } \left(\boldsymbol{x}_{k-q}-\boldsymbol{x}_{k-q-1}\right)\text{,} \dots\text{, }\left(\boldsymbol{x}_{k}-\boldsymbol{x}_{k-1}\right)\right] \text{,} \\
	R_{k} &= \left[\left(\boldsymbol{r}_{k-q+1}-\boldsymbol{r}_{k-q}\right)\text{, } \left(\boldsymbol{r}_{k-q}-\boldsymbol{r}_{k-q-1}\right)\text{,} \dots\text{, }\left(\boldsymbol{r}_{k}-\boldsymbol{r}_{k-1}\right)\right]\text{.} 
\end{align*}
Rewriting Eq. \eqref{Eq3.3} as
$$ \overline{\boldsymbol{x}}_{k} =\boldsymbol{x} _{k}- X_{k} \boldsymbol{\Gamma}_{q}\text{,}   $$
the corresponding modified residual is $ \overline {\boldsymbol{r} }_{k} =  M^{-1} \left(\boldsymbol{b}-A \overline{ \boldsymbol{x}}_{k} \right) $. 
The parameter \textcolor{black}{vector} $ \boldsymbol{\Gamma}_{q} $ \textcolor{black}{is} chosen by minimizing the $ l_{2} $ norm of $ \overline {\boldsymbol{r}} _{k} $ as follows:
\begin{align*}
	\boldsymbol{\Gamma}_{q} &= \arg\min  \Vert M^{-1} \left(\boldsymbol{b}-A \overline {\boldsymbol{x}}_{k} \right)  \Vert _{2} \\
	&= \arg\min  \Vert \boldsymbol{r}_{k} -R_{k} \boldsymbol{\Gamma}_{q} \Vert _{2}\text{.}
\end{align*} 
\textcolor{black}{If} the columns of $ R_{k} $ are linearly independent, the optimal solution $ \boldsymbol{\Gamma}_{q} $ satisfies
$$  \boldsymbol{\Gamma}_{q} = \left(R_{k}^{\top}R_{k}\right)^{-1}R_{k}^{\top} \boldsymbol{r}_{k}\text{.} $$

The second step is to generalize Eq. \eqref{Eq3.2} using $ \overline{ \boldsymbol{x}}_{k}\text{, }  \overline{\boldsymbol{r}}_{k}  $.
The AR method is
$$ \boldsymbol{x}_{k+1} =\overline{\boldsymbol{x}} _{k}+\beta_{k} \overline{\boldsymbol{r}} _{k}\text{,} $$
which is written as
$$  \boldsymbol{x}_{k+1}=\boldsymbol{x}_{k}+\left[\beta_{k} I-\left(X_{k}+\beta_{k} R_{k}\right)\left(R_{k}^{\top}R_{k}\right)^{-1}R_{k}^{\top}\right] \boldsymbol{r}_{k}\text{.} $$

The relaxation parameter $ \beta_{k}  $ is chosen by minimizing the quasi\textcolor{black}{-}residual as follows:
\begin{align*}
	\beta_{k} &= \arg\min  \Vert  \boldsymbol{b}-A \boldsymbol{x}_{k+1}  \Vert _{2} \\
	&= \arg\min  \Vert \boldsymbol{b}-A\overline{\boldsymbol{x}}_{k} -\beta_{k} A \overline{\boldsymbol{r}}_{k} \Vert _{2}\text{.}
\end{align*} 
The optimal solution $ \beta_{k} $ satisfies
$$ \beta_{k} = \frac{\left( \boldsymbol{b}-A\overline{\boldsymbol{x}}_{k}\text{, } A \overline{\boldsymbol{r}}_{k}\right)}{\left(A\overline{\boldsymbol{r}}_{k}\text{, } A\overline{\boldsymbol{r}}_{k}\right)} \text{.}$$

Denoting $ B_{k}=\beta_{k} I-\left(X_{k}+\beta_{k} R_{k}\right)\left(R_{k}^{\top}R_{k}\right)^{-1}R_{k}^{\top} $, the error of the AR iteration is $ \boldsymbol{e}_{k+1} = \left(I-B_{k} M^{-1}A\right) \boldsymbol{e}_{k} $, 
the residual is $ \boldsymbol{r}_{k+1} = \left(I-B_{k} M^{-1}A\right)\boldsymbol{r}_{k} $, and the convergence requires $$ \rho \left( I-B_{k} M^{-1}A \right) <1 \text{.} $$
It can be concluded that the AR iteration converges faster if $ B_{k} $ is better approximated to $ A^{-1}M $.

\subsubsection{Alternating Anderson--Richardson method}
``Low frequency'' and ``high frequency'' residual components, respectively, represent the eigenvalues of \textcolor{black}{$ (I-\omega_{k}$$ M^{-1}A) $} with values close to unity and zero.
The weighted Jacobi method cannot efficiently reduce the low-frequency components.
To fix this problem,
Pratapa et al. \cite{r16} proposed the AAJ method and generalized it to the AAR method in \cite{r17}.
The AAR method incorporates the Richardson and AR methods, the AAR method is
$$ \boldsymbol{x}_{k+1}=\boldsymbol{x}_{k}+V_{k} \boldsymbol{r}_{k}\text{,} \quad  k=0\text{, } 1\text{,} \dots \text{,} $$
where $ V _{k} $ is
\begin{equation*} 
	V_{k}=\left\{
	\begin{aligned}
		& \omega_{k} I & \text{if} \ \frac{k+1}{p} \notin \mathbb{N}\text{,} \\
		& \beta_{k} I-\left(X_{k}+\beta_{k} R_{k}\right)\left(R_{k}^{\top}R_{k}\right)^{-1}R_{k}^{\top}  & \text{if} \ \frac{k+1}{p} \in \mathbb{N}\text{.}\
	\end{aligned}
	\right.
\end{equation*}

The pseudocode of the AAR method for solving linear systems is presented in Algorithm \ref{alg1}.
\begin{algorithm}
	\renewcommand{\algorithmicrequire}{\textbf{Input:}}
	\renewcommand{\algorithmicensure}{\textbf{Output:}}
	\caption{The AAR method for linear systems}
	\label{alg1}
	\begin{algorithmic}[1]
		\REQUIRE $A\text{, } M\text{, } \boldsymbol{b}\text{, } \boldsymbol{x}_{0}\text{, } p\text{, } q\text{, } k_{max}\text{, } tol$;
		\ENSURE  $ \boldsymbol{x}_{k+1}$.
		\STATE initial $ k=0 $, $\boldsymbol{x}_{old} = \boldsymbol{x}_{0}$;
		\WHILE{$ k \leqslant k_{max} $ and $ \frac{\Vert \boldsymbol{b} - A \boldsymbol{x}_{k} \Vert _2}{\Vert \boldsymbol{b} - A \boldsymbol{x}_{0} \Vert _2} > tol $}
		\STATE $  \boldsymbol{r}_{k} =M^{-1}\left(\boldsymbol{b} -A \boldsymbol{x}_{k}^{m-1} \right)   $;
		\IF{$ k >1 $}
		\STATE $\boldsymbol{x}\left(:\text{, }\text{mod}\left(k-2\text{, }q\right)+1\right) =\boldsymbol{x}_{k} - \boldsymbol{x}_{old} $; (``mod" is the remainder operation)
		\STATE $R\left(:\text{, }\text{mod}\left(k-2\text{, }q\right)+1\right) = \boldsymbol{r}_{k} - \boldsymbol{r}_{old} $;
		\ENDIF
		\STATE $\boldsymbol{x}_{old} = \boldsymbol{x}_{k}\text{; } \boldsymbol{r}_{old} = \boldsymbol{r}_{k} $;
		\IF{$\frac{k+1}{p} \notin \mathbb{N}$}
		\STATE$ \omega_{k} = \frac{\left(\boldsymbol{b}-A\boldsymbol{x}_{k}\text{, } A\boldsymbol{r}_{k}\right)}{\left(A\boldsymbol{r}_{k}\text{, } A\boldsymbol{r}_{k}\right)} $; $\boldsymbol{x}_{k+1} =\boldsymbol{x}_{k}+\omega_{k} \boldsymbol{r}_{k} $;
		\ELSE
		\STATE $ \boldsymbol{\Gamma}_{k} = \left(R_{k}^{\top}R_{k}\right)^{-1}R_{k}^{\top} \boldsymbol{r}_{k} $;
		\STATE $ \overline{\boldsymbol{x}} _{k} =  \boldsymbol{x} _{k} - X_{k} \boldsymbol{\Gamma}_{k} $; $ \overline{\boldsymbol{r}} _{k} =M^{-1}\left(\boldsymbol{b}-A\overline{\boldsymbol{x}} _{k}\right) $;		
		\STATE $ \beta_{k} = \frac{\left( \boldsymbol{b}-A\overline{\boldsymbol{x}}_{k}\text{, } A \overline{\boldsymbol{r}}_{k}\right)}{\left(A\overline{ \boldsymbol{r}}_{k} \text{, } A\overline{ \boldsymbol{r}}_{k}\right)} $; $ \boldsymbol{x}_{k+1} =\overline{\boldsymbol{x}} _{k}+\beta_{k} \overline{\boldsymbol{r}} _{k} $;
		\ENDIF
		\STATE $ \boldsymbol{x}_{old} = \boldsymbol{x}_{k} ; \boldsymbol{x}_{k} = \boldsymbol{x}_{k+1} $;
		\ENDWHILE
	\end{algorithmic}  
\end{algorithm}

\subsection{Tensor Alternating Anderson--Richardson method}
\label{subsec:subsec3.2}

We apply the idea of the AAR method for solving Eq. \eqref{Eq1.1} and propose a TAAR method.
We first present a TR method based on tensor regular splittings in Section \ref{subsubsec:subsubsec3.2.1}. Then, we apply Anderson acceleration to the TR method and derive a TAR method in Section \ref{subsubsec:subsubsec3.2.2}.
Finally, we conclude the derivation of the TAAR method in Section \ref{subsubsec:subsubsec3.2.3}.
We compare the computational cost between the TAAR method and three tensor splitting methods in Section \ref{subsubsec:subsubsec3.2.4}.

\subsubsection{Tensor Richardson method}
\label{subsubsec:subsubsec3.2.1}

Considering Eq. \eqref{Eq1.1}
with $ \boldsymbol{b} \in \mathbb{R}_{++}^{n} $ and a nonsingular $ \mathcal{M} $-tensor $\mathcal{A} $, 
$\mathcal{A} $ has a regular splitting
$$ \mathcal{A}=\mathcal{E}-\mathcal{F} \text{,}$$
where \textcolor{black}{a} left-nonsingular tensor $ \mathcal{E} \in \mathbb{R}^{[m\text{,}n]}  $ satisfies $ \mathcal{E}=\mathit{M}\left(\mathcal{E}\right) \mathcal{I}  $ with  $  \mathit{M}\left(\mathcal{E}\right)_{ij} = \mathcal{E} _{ij \dots j}\text{, }  i\text{, }j=1\text{, } 2\text{,} \dots \text{, } n $. 
Using this regular splitting, we rewrite Eq. \eqref{Eq1.1} as follows:
$$ \boldsymbol{x}^{\left[m-1\right]} = \mathit{M}\left(\mathcal{E}\right) ^{-1} \left(\boldsymbol{b}+\mathcal{F} \boldsymbol{x}^{m-1} \right)=:g\left(\boldsymbol{x}\right) \text{.}$$
The residual is defined as
\begin{align*}
	r \left(\boldsymbol{x}\right)  &= g\left(\boldsymbol{x}\right)-\boldsymbol{x}^{\left[m-1\right]}\\
	&= g\left(\boldsymbol{x}\right)- \mathit{M}\left(\mathcal{E}\right) ^{-1} \mathcal{E} \boldsymbol{x}^{m-1}\\
	&= \mathit{M}\left(\mathcal{E}\right) ^{-1} \left(\boldsymbol{b}-\mathcal{A}\boldsymbol{x}^{m-1} \right) \text{.}
\end{align*}

The TR method is
\begin{equation}
	\label{E3.2.1}
	\boldsymbol{x}_{k+1}^{\left[m-1\right]} = \boldsymbol{x}_{k}^{\left[m-1\right]} + \omega _{k} \boldsymbol{r}_{k}
\end{equation}
with the residual $ \boldsymbol{r}_{k}= \mathit{M}\left(\mathcal{E}\right) ^{-1} \left(\boldsymbol{b}-\mathcal{A}\boldsymbol{x}_{k}^{m-1} \right) $ and the relaxation parameter $ \omega _{k} \in \mathbb{R} $.
There are three preconditioners $ \mathit{M}\left(\mathcal{E}\right) $ according to \cite{r6}:
\begin{enumerate}[(1)]
	\item Jacobi-type preconditioner (PJ): $ \mathit{M}\left(\mathcal{E}\right)= D\left( \mathit{M}\left(\mathcal{A}\right) \right) $, where $ D\left( \mathit{M}\left(\mathcal{A}\right) \right) $ is the diagonal part of $ \mathit{M}\left(\mathcal{A}\right) $; 
	\item Gauss--Seidel-type preconditioner (PGS): $ \mathit{M}\left(\mathcal{E}\right)= L\left( \mathit{M}\left(\mathcal{A}\right)\right) $, where $ L\left( \mathit{M}\left(\mathcal{A}\right)\right) $ is the lower triangular part of $ \mathit{M}\left(\mathcal{A}\right) $;
	\item FULLM-type preconditioner (PF): $ \mathit{M}\left(\mathcal{E}\right)= \mathit{M}\left(\mathcal{A}\right) $.
\end{enumerate}

We can obtain the formula of $ \boldsymbol{x}_{k+1} $ from Eq. \eqref{E3.2.1}.
However, if we choose the relaxation parameter $ \omega_{k} $ by minimizing $ \boldsymbol{r}_{k+1} $ directly,
$$ \omega_{k} =  \arg\min  \Vert \mathit{M}\left(\mathcal{E}\right) ^{-1} \left(\boldsymbol{b}-\mathcal{A}\boldsymbol{x}_{k+1}^{m-1} \right)  \Vert _{2} \text{,} $$
it is inconvenient to solve the above minimal problem because the product operation between a tensor and vector is nonlinear.
We need to find another technique to determine $ \omega_{k} $.

Let $ B = \mathcal{A} \boldsymbol{x}_{k}^{m-2} $, we approximate $ \mathcal{A}\boldsymbol{x}_{k+1}^{m-1} $ with $B  \boldsymbol{x}_{k+1} $ 
and choose $ \omega_{k} $ by minimizing the approximate quasi\textcolor{black}{-}residual
\begin{equation} 
	\begin{aligned}
		\label{Eq3.2.2}
		\omega_{k} = \arg\min  \Vert \boldsymbol{b} -B \boldsymbol{x}_{k+1} \Vert _{2}.
	\end{aligned} 
\end{equation}  
Eq. \eqref{Eq3.2.2} needs to calculate $ \boldsymbol{x}_{k+1} $.  We reformulate Eq. \eqref{E3.2.1} as follows:
$$ \boldsymbol{x}_{k+1}^{[m-1]} - \boldsymbol{x}_{k}^{[m-1]} = \omega _{k} \boldsymbol{r}_{k}\text{. } $$
The left-hand side is
\begin{equation}
	\label{Eq3.2.2s}
	\boldsymbol{x}_{k+1}^{\left[m-1\right]} - \boldsymbol{x}_{k}^{\left[m-1\right]}  = \left(\boldsymbol{x}_{k+1}-\boldsymbol{x}_{k} \right).*\left( \boldsymbol{x}_{k+1}^{\left[m-2\right]} +  \boldsymbol{x}_{k+1}^{\left[m-3\right]} .* \boldsymbol{x}_{k} + \dots + \boldsymbol{x}_{k}^{\left[m-2\right]}  \right) \text{,}
\end{equation}
where $ \boldsymbol{x}_{k}^{[m-j]} = \left[\left(\boldsymbol{x}_{k}\right)_{1}^{m-j}\text{,}\left(\boldsymbol{x}_{k}\right)_{2}^{m-j}\text{,}\dots\text{,}\left(\boldsymbol{x}_{k}\right)_{n}^{m-j}\right]^{\top}\text{, } j=1\text{, } 2\text{,} \dots\text{, } m-1 $.
For any $ \boldsymbol{x} $, $ \boldsymbol{y} \in \mathbb{R}^{n}$, $ \boldsymbol{x} .* \boldsymbol{y} =  \left[   x_{1}y_{1} \text{, } x_{2}y_{2} \text{,} \dots \text{, } x_{n}y_{n}\right] ^{\textcolor{black}{\top}}$.
We approximate $\left( \boldsymbol{x}_{k+1}^{\left[m-2\right]} +  \boldsymbol{x}_{k+1}^{\left[m-3\right]}\textcolor{black}{.* } \boldsymbol{x}_{k} + \dots + \boldsymbol{x}_{k}^{\left[m-2\right]}  \right)$ by 
substituting $ \boldsymbol{x}_{k+1}^{[m-j]} $ with $ \boldsymbol{x}_{k}^{[m-j]} $,
\begin{equation}  \nonumber
	\boldsymbol{x}_{k+1}^{\left[m-1\right]}-\boldsymbol{x}_{k}^{\left[m-1\right]} \approx (m-1)\left(\boldsymbol{x}_{k+1}-\boldsymbol{x}_{k}\right).*\boldsymbol{x}_{k}^{\left[m-2\right]}\text{.}
\end{equation}
Introducing a parameter $  \tilde{\eta}_{k} $ and denoting $ \eta_{k} = (m-1)  \tilde{\eta}_{k}  $, the above equation becomes 
\begin{equation} 
	\label{Eq3.2.3.1}
	\boldsymbol{x}_{k+1}^{\left[m-1\right]}-\boldsymbol{x}_{k}^{\left[m-1\right]} \approx \eta_{k} \left(\boldsymbol{x}_{k+1}-\boldsymbol{x}_{k}\right).*\boldsymbol{x}_{k}^{\left[m-2\right]}\text{. }
\end{equation} 
According to Eqs. \eqref{E3.2.1} and \eqref{Eq3.2.3.1}, $ x_{k+1} $ is
$$ 	\boldsymbol{x}_{k+1} = \boldsymbol{x}_{k} + \frac{\omega _{k}}{\eta_{k}} \boldsymbol{u}_{1}\text{,} $$
where $  \boldsymbol{u}_{1} = \boldsymbol{r}_{k}./\boldsymbol{x}_{k}^{\left[m-2\right]} $ with $ \boldsymbol{r}_{k}./\boldsymbol{x}_{k}^{\left[m-2\right]} = \left[ \frac{\left(\boldsymbol{r}_{k}\right)_{1}}{\left(\boldsymbol{x}_{k}\right)_{1}^{[m-2]}}\text{, } \frac{\left(\boldsymbol{r}_{k}\right)_{2}}{\left(\boldsymbol{x}_{k}\right)_{2}^{[m-2]}}\text{,} \dots\text{, } \frac{\left(\boldsymbol{r}_{k}\right)_{n}}{\left(\boldsymbol{x}_{k}\right)_{n}^{[m-2]}}   \right]^{\top}$.
Eq. \eqref{Eq3.2.2} is rewritten as
\begin{equation*}  
	\begin{aligned}
		\frac{\omega _{k}}{\eta_{k}}  &=\arg\min  \Vert \boldsymbol{b} -B \left(  \boldsymbol{x}_{k} + \frac{\omega _{k}}{\eta_{k}}  \boldsymbol{u}_{1} \right) \Vert _{2} \\
		&= \arg\min  \Vert \boldsymbol{b} - \mathcal{A} \boldsymbol{x}_{k}^{m-1} -\frac{\omega _{k}}{\eta_{k}}  B  \boldsymbol{u}_{1} \Vert _{2}\text{. } 
	\end{aligned} 
\end{equation*} 
We want to choose $\omega_{k}$ by minimzing the approximate quasi\textcolor{black}{-}residual.\vspace{0.2ex}
The parameter value that minimizes the approximate quasi\textcolor{black}{-}residual is denoted as $\frac{\omega _{k}}{\eta_{k}}  $.
Therefore, we take the value of $ \frac{\omega _{k}}{\eta_{k}}  $ as $ \omega_{k} $, which means $ \eta _{k} =1 $.
Eq. \eqref{Eq3.2.3.1} becomes
\begin{equation} 
	\label{Eq3.2.3}
	\boldsymbol{x}_{k+1}^{\left[m-1\right]}-\boldsymbol{x}_{k}^{\left[m-1\right]} \approx \left(\boldsymbol{x}_{k+1}-\boldsymbol{x}_{k}\right).*\boldsymbol{x}_{k}^{\left[m-2\right]}\text{. }
\end{equation}
We can obtain $ \omega_{k} $ as
$$ \omega_{k} = \frac{\left(\boldsymbol{b} - \mathcal{A} \boldsymbol{x}_{k}^{m-1}\text{, }B \boldsymbol{u}_{1}\right)}{\left(B  \boldsymbol{u}_{1}\text{, }B  \boldsymbol{u}_{1}\right)} \text{.}$$

\subsubsection{Tensor Anderson--Richardson method}
\label{subsubsec:subsubsec3.2.2}
Now, we apply Anderson acceleration to the TR method for solving multilinear systems and derive a TAR method. The TAR method includes two steps:

The first step is to correct $ \boldsymbol{x}_{k}^{\left[m-1\right]}$ and $ \boldsymbol{r}_{k} $ in Eq. \eqref{E3.2.1}. 
In Section \ref{subsubsec:subsubsec3.1.2}, Pratapa et al. \cite{r17} used \textcolor{black} {the parameter vector} $ \boldsymbol{\Gamma}_{q} \in \mathbb{R}^{q} $ to correct $ \boldsymbol{x}_{k} $.
Considering the current iteration number $ k $ may be smaller than $ q $ when running the TAR method, and let $ l=\min \{ q\text{, } k\} $, 
we use parameters $\gamma_{1}\text{, } \gamma_{2}\text{,} \dots\text{, } \gamma_{l} $ to correct the $ \boldsymbol{x}_{k}^{\left[m-1\right]} $ in Eq. \eqref{E3.2.1}, denoted by $ \overline{\boldsymbol{x}} _{k} $,
\begin{equation}
	\label{Eq3.2.4}
	\overline{\boldsymbol{x}}_{k} =\boldsymbol{x}_{k}^{\left[m-1\right]} - \sum_{j=1}^{l} \gamma_{j} \left( \boldsymbol{x}_{k-l+j} - \boldsymbol{x}_{k-l+j-1}\right)\text{.}
\end{equation} 
Let
\begin{align*}
	X_{k} &=\left[\left(\boldsymbol{x}_{k-l+1}^{\left[m-1\right]}-\boldsymbol{x}_{k-l}^{\left[m-1\right]}\right)\text{,} \left(\boldsymbol{x}_{k-l}^{\left[m-1\right]}-\boldsymbol{x}_{k-l-1}^{\left[m-1\right]}\right)\text{,} \dots\text{,} \left(\boldsymbol{x}_{k+1}^{\left[m-1\right]}-\boldsymbol{x}_{k}^{\left[m-1\right]}\right)\right] \in \mathbb{R}^{n \times l}\text{,}  \\
	R_{k} & = \left[\left(\boldsymbol{r}_{k-l+1}-\boldsymbol{r}_{k-l}\right)\text{,} \left(\boldsymbol{r}_{k-l}-\boldsymbol{r}_{k-l-1}\right)\text{,} \dots\text{,}\left(\boldsymbol{r}_{k}-\boldsymbol{r}_{k-1}\right)\right]\in \mathbb{R}^{n \times l}\text{,}  \\
	\boldsymbol{\Gamma}_{l} &= \left[\gamma_{1}\text{, } \gamma_{2}\text{,} \dots\text{, } \gamma_{l} \right]^{\top} \in \mathbb{R}^{l} \text{,}
\end{align*}
We rewrite Eq. \eqref{Eq3.2.4} as 
$$	\overline{\boldsymbol{x}} _{k} = \boldsymbol{x}_{k}^{\left[m-1\right]} - X_{k} \boldsymbol{\Gamma}_{l} \text{.}$$

If we choose $ \boldsymbol{\Gamma}_{l} $ by directly minimizing $ \mathit{M}\left(\mathcal{E}\right) ^{-1} \left(\boldsymbol{b}-\mathcal{A}  \overline{\boldsymbol{x}}_{k}^{m-1} \right) $, 
it is difficult to obtain the optimal solution because the product operation between tensor and vector is nonlinear.
Recalling that the modified residual in a linear system is 
\begin{equation}  \nonumber
	\begin{aligned}
		\overline{\boldsymbol{r}}_{k}  &=\ M^{-1} \left(\boldsymbol{b}-A \overline{\boldsymbol{x}}_{k}\right) \\
		&= \boldsymbol{r}_{k}-R_{k}\boldsymbol{\Gamma}_{q} \text{, } 
	\end{aligned} 
\end{equation}  
we define the modified residual of a multilinear system as
$ \overline{\boldsymbol{r}}_{k} =\boldsymbol{r}_{k} - R_{k} \boldsymbol{\Gamma}_{l} $ and choose $ \boldsymbol{\Gamma}_{l} $ by minimizing $ \overline{\boldsymbol{r}}_{k}  $,
$$  \boldsymbol{\Gamma}_{l} =\arg \min \Vert \boldsymbol{r}_{k}-R_{k} \boldsymbol{\Gamma}_{l}  \Vert_{2} \text{.} $$
\textcolor{black}{If} the columns of $ R_{k} $ are linearly independent, the optimal solution $ \boldsymbol{\Gamma}_{l} $ satisfies
$$  \boldsymbol{\Gamma}_{l} = \left(R_{k}^{\top}R_{k}\right)^{-1}R_{k}^{\top} \boldsymbol{r}_{k}\text{.} $$

The second step is to generalize Eq. \eqref{E3.2.1} using $ \overline{\boldsymbol{x}}_{k}\text{ and }  \overline{\boldsymbol{r}}_{k}  $,
$$ \boldsymbol{x}_{k+1}^{\left[m-1\right]} =\overline{\boldsymbol{x}} _{k}+\beta_{k} \overline{\boldsymbol{r}} _{k}\text{,} $$
then the TAR method is
\begin{equation} 
	\label{Eq3.2.5}
	\boldsymbol{x}_{k+1}^{\left[m-1\right]}=\boldsymbol{x}_{k}^{[m-1]}+\left[\beta_{k} I-\left(X_{k}+\beta_{k}  R_{k}\right)\left(R_{k}^{\top}R_{k}\right)^{-1}R_{k}^{\top}\right] \boldsymbol{r}_{k} \text{.}
\end{equation}

We choose $ \beta_{k} $ in the same manner as choosing $ \omega_{k} $ in Section \ref{subsubsec:subsubsec3.2.1}. According to Eqs. \eqref{Eq3.2.3} and \eqref{Eq3.2.5}, we obtain the following formula
\begin{equation} \nonumber
	\label{e5.6}
	\boldsymbol{x}_{k+1} =\boldsymbol{x}_{k}-\boldsymbol{u}_{2}+\beta_{k} \boldsymbol{u}_{3}\text{,}
\end{equation}
where 
\begin{align*}
	\boldsymbol{u}_{2} &=\left(X_{k} \boldsymbol{\Gamma}_{l}\right)./ \boldsymbol{x}_{k}^{\left[m-2\right]} \text{,} \\
	\boldsymbol{u}_{3} &=\left(\boldsymbol{r}_{k}-R_{k}\boldsymbol{\Gamma}_{l}\right)./ \boldsymbol{x}_{k}^{\left[m-2\right]} 
\end{align*}
with 
\begin{align*}
	\left(X_{k} \boldsymbol{\Gamma}_{l}\right)./\boldsymbol{x}_{k}^{\left[m-2\right]} &= \left[ \frac{\left(X_{k} \boldsymbol{\Gamma}_{l}\right)_{1}}{\left(\boldsymbol{x}_{k}\right)_{1} ^{\left[m-2\right]}}\text{, }  \frac{\left(X_{k}\boldsymbol{\Gamma}_{l}\right)_{2}}{\left(\boldsymbol{x}_{k}\right)_{2} ^{\left[m-2\right]}}\text{,} \dots\text{, } \frac{\left(X_{k} \boldsymbol{\Gamma}_{l}\right)_{n}}{\left(\boldsymbol{x}_{k}\right)_{n} ^{\left[m-2\right]} }  \right]^{\top} \text{,}\\
	\left(\boldsymbol{r}_{k}-R_{k}\boldsymbol{\Gamma}_{l}\right)./ \boldsymbol{x}_{k}^{[m-2]} &= \left[ \frac{\left(\boldsymbol{r}_{k}\right)_{1}-\left(R_{k}\boldsymbol{\Gamma}_{l}\right)_{1}}{\left(\boldsymbol{x}_{k}\right)^{\left[m-2\right]}_{1}}\text{, } \frac{\left(\boldsymbol{r}_{k}\right)_{2}-\left(R_{k}\boldsymbol{\Gamma}_{l}\right)_{2}}{\left(\boldsymbol{x}_{k}\right)^{\left[m-2\right]}_{2}}\text{,} \dots\text{, } \frac{\left(\boldsymbol{r}_{k}\right)_{n}-\left(R_{k}\boldsymbol{\Gamma}_{l}\right)_{n}}{\left(\boldsymbol{x}_{k}\right) ^{\left[m-2\right]}_{n}}   \right]^{\top}\text{.}
\end{align*}
We choose $ \beta_{k} $ by minimizing \textcolor{black}{the} approximate quasi\textcolor{black}{-}residual as follows:
\begin{equation}  \nonumber
	\begin{aligned}
		\beta_{k} &=\arg\min  \Vert \boldsymbol{b} -B \boldsymbol{x}_{k+1}  \Vert _{2} \\
		&=\arg\min  \Vert \boldsymbol{b} -B \left( \boldsymbol{x}_{k}-\boldsymbol{u}_{2}+\beta_{k} \boldsymbol{u}_{3} \right) \Vert _{2} \text{,} 
	\end{aligned} 
\end{equation} 
where $ B = \mathcal{A} \boldsymbol{x}_{k}^{m-2} $ is a matrix.  
The optimal solution $ \beta_{k}  $ satisfies
$$ \beta_{k} = \frac{\left(\boldsymbol{b} - \mathcal{A} \boldsymbol{x}_{k}^{m-1}+ B \boldsymbol{u}_{2}\text{, } B  \boldsymbol{u}_{3}\right)}{\left(B  \boldsymbol{u}_{3}\text{, } B  \boldsymbol{u}_{3}\right)} \text{.}$$

\subsubsection{Tensor Alternating Anderson--Richardson method}
\label{subsubsec:subsubsec3.2.3}

We apply one iteration of TAR after every $\left(p-1\right)$ iterations of TR and propose a TAAR method. The TAAR method is
\begin{equation} 
	\boldsymbol{x}_{k+1}^{[m-1]}=\boldsymbol{x}_{k}^{\left[m-1\right]}+V_{k}\boldsymbol{r}_{k}\text{,}\quad  k=0\text{, } 1\text{,} \dots\text{,}
\end{equation}
where residual ${r}_{k}$ satisfies $ \boldsymbol{r}_{k}= \mathit{M}\left(\mathcal{E}\right) ^{-1} \left(\boldsymbol{b}-\mathcal{A}\boldsymbol{x}_{k}^{m-1} \right) $ and $ V _{k} $ is
\begin{equation} \nonumber
	V_{k}=\left\{
	\begin{aligned}
		& \omega_{k} I & \text{if} \ \frac{k+1}{p} \notin \mathbb{N}\text{,} \\
		& \beta_{k} I-\left(X_{k}+\beta_{k} R_{k}\right)\left(R_{k}^{\top}R_{k}\right)^{-1}R_{k}^{\top}  & \text{if} \ \frac{k+1}{p} \in \mathbb{N}\text{.}
	\end{aligned}
	\right.
\end{equation}

In the TAAR method, we first need to produce the preconditioner $ \mathit{M}\left(\mathcal{E}\right) $, elementwise,
$$  \mathit{M}\left(\mathcal{E}\right)_{ij} = \mathcal{E} _{ij \dots j}\text{,} \quad  i\text{, }j=1\text{, }2\text{,} \dots \text{, } n \text{.}$$
We use the tensor toolbox 3.2.1 \cite{r27} and produce the tensor index $ \left(i\text{, } j\text{,} \dots\text{, } j \right) $ using the matrix index $ \left(i\text{, } j \right) $. 
The specific algorithm is given by Algorithm \ref{alg2}.
\begin{algorithm}
	\renewcommand{\algorithmicrequire}{\textbf{Input:}}
	\renewcommand{\algorithmicensure}{\textbf{Output:}}
	\caption{Produce $ \mathit{M}\left(\mathcal{E}\right) $}
	\label{alg2}
	\begin{algorithmic}[1]
		\REQUIRE $\mathcal{A}\text{, } m\text{, } n$;
		\ENSURE  $ \mathit{M}\left(\mathcal{E}\right) $.
		\STATE $ msize $ is the size of $ \mathit{M}\left(\mathcal{A}\right) $;
		\STATE $ \left[i\text{, } j\right] = \text{ind2sub}\left(msize\text{, } 1:n^{2}\right) $; (``ind2sub" returns the index of  $\mathit{M}\left(\mathcal{A}\right)$) 
		\STATE produce $ index= \left[i\text{, } j\text{,} \dots\text{, } j\right] $;
		\STATE  construct $ \mathit{M}\left(\mathcal{A}\right)_{ij} = a _{ij \dots j}  $;
		\IF{choose the Jacobi-type preconditioner}
		\STATE $ \mathit{M}\left(\mathcal{E}\right)= D\left( \mathit{M}\left(\mathcal{A}\right)\right) $; (``$D\left( \mathit{M}\left(\mathcal{A}\right)\right)$" denotes the diagonal part of $ \mathit{M}\left(\mathcal{A}\right) $)
		\ENDIF
		\IF{choose the Gauss-Seidel-type preconditioner}
		\STATE $ \mathit{M}\left(\mathcal{E}\right)= \text{tril}\left( \mathit{M}\left(\mathcal{A}\right)\right) $; (``tril" returns the lower triangular part of $ \mathit{M}\left(\mathcal{A}\right) $)
		\ENDIF
		\IF{choose the FULLM-type preconditioner}
		\STATE $ \mathit{M}\left(\mathcal{E}\right)= \mathit{M}\left(\mathcal{A}\right)  $;
		\ENDIF  
	\end{algorithmic}  
\end{algorithm}

The pseudocode of the TAAR method for multilinear systems is given \textcolor{black}{in} Algorithm \ref{alg3}.
\begin{algorithm}
	\renewcommand{\algorithmicrequire}{\textbf{Input:}}
	\renewcommand{\algorithmicensure}{\textbf{Output:}}
	\caption{The TAAR method for multilinear systems}
	\label{alg3}
	\begin{algorithmic}[H]
		\REQUIRE $ \mathcal{A}\text{, } \boldsymbol{b}\text{, } m\text{, } n\text{, } \boldsymbol{x}_{0}\text{, } p\text{, } q\text{, } k_{max}\text{, } tol$;
		\ENSURE  $ \boldsymbol{x}_{k+1}$.
		\STATE initial $ k=0 $, $\boldsymbol{x}_{old} =\boldsymbol{x}_{0}$; 
		\STATE produce $ \mathit{M}\left(\mathcal{E}\right) $ according to Algorithm \ref{alg2}; 
		\WHILE{$ k \leqslant k_{max} $ and $ \frac{\Vert \boldsymbol{b} - \mathcal{A} \boldsymbol{x}_{k} ^{m-1}\Vert _2}{\Vert \boldsymbol{b} - \mathcal{A} \boldsymbol{x}_{0} ^{m-1} \Vert _2} > tol $}
		\STATE $  \boldsymbol{r}_{k} =\mathit{M}\left(\mathcal{E}\right) ^{-1} \left(\boldsymbol{b}-\mathcal{A}\boldsymbol{x}_{k}^{m-1} \right) $;
		\IF{$ k >1 $}
		\STATE $ l=\min \{ q\text{, }k\} \text{, }$
		\STATE $X\left(:\text{, }\text{mod}\left(k-2\text{, }q\right)+1\right) =\boldsymbol{x}_{k} ^{[m-1]}- \boldsymbol{x}_{old}^{[m-1]} $; 
		\STATE $R\left(:\text{, }\text{mod}\left(k-2\text{, }q\right)+1\right) = \boldsymbol{r}_{k} - \boldsymbol{r}_{old}$;
		\ENDIF
		\STATE $\boldsymbol{x}_{old}^{\left[m-1\right]} = \boldsymbol{x}_{k}^{\left[m-1\right]}; \boldsymbol{r}_{old} = \boldsymbol{r}_{k} $;
		\IF{$\frac{k+1}{p} \notin \mathbb{N}$}
		\STATE $  \boldsymbol{u}_{1} = \boldsymbol{r}_{k}./\boldsymbol{x}_{k}^{\left[m-2\right]} \text{; } \omega_{k} = \frac{\left(\boldsymbol{b} - \mathcal{A} \boldsymbol{x}_{k}^{m-1}\text{, }B  u_{1}\right)}{\left(B  u_{1}\text{, }B  u_{1}\right)} $;
		\STATE $\boldsymbol{x}_{k+1} ^{[m-1]} =\boldsymbol{x}_{k}^{[m-1]}+\omega_{k} \boldsymbol{r}_{k} $;
		\ELSE
		\STATE $ \boldsymbol{\Gamma}_{l} = \left(R_{k}^{\top}R_{k}\right)^{-1}R_{k}^{\top} \boldsymbol{r}_{k} $, $\boldsymbol{u}_{2}=\left(X_{k} \boldsymbol{\Gamma}_{l}\right)./ \boldsymbol{x}_{k}^{\left[m-2\right]} \text{; } \boldsymbol{u}_{3}=\left(\boldsymbol{r}_{k}-R_{k}\boldsymbol{\Gamma}_{l}\right)./ \boldsymbol{x}_{k}^{[m-2]} $;
		\STATE $  \beta_{k} = \frac{\left(\boldsymbol{b} - \mathcal{A} \boldsymbol{x}_{k}^{m-1}+ B \boldsymbol{u}_{2}\text{, }B \boldsymbol{u}_{3}\right)}{\left(B  \boldsymbol{u}_{3}\text{, }B  \boldsymbol{u}_{3}\right)}   $;
		\STATE $ \boldsymbol{x}_{k+1}^{\left[m-1\right]}=\boldsymbol{x}_{k}^{[m-1]}+[\beta_{k} I-\left(X_{k}+\beta_{k}  R_{k}\right)\left(R_{k}^{\top}R_{k}\right)^{-1}R_{k}^{\top}] \boldsymbol{r}_{k} $;		
		\ENDIF
		\STATE $ \boldsymbol{x}_{old} ^{\left[m-1\right]}= \boldsymbol{x}_{k} ^{\left[m-1\right]} \text{; } \boldsymbol{x}_{k}^{\left[m-1\right]} = \boldsymbol{x}_{k+1}^{\left[m-1\right]} $;
		\ENDWHILE
	\end{algorithmic}  
\end{algorithm}

\subsubsection{Computational cost}
\label{subsubsec:subsubsec3.2.4}
The computational cost corresponds to one operation on two floating point numbers.
The computational cost of each iteration in the TAAR method has two components, the computational cost of TR and the computational cost of TAR.
The TR and TAR methods need to compute $ \mathit{M}\left(\mathcal{E}\right) ^{-1} \left(\boldsymbol{b}-\mathcal{A}\boldsymbol{x}_{k}^{m-1}\right) $.
The computational \textcolor{black}{cost} of \textcolor{black}{$ \mathit{M}\left(\mathcal{E}\right) ^{-1} (\boldsymbol{b}-\mathcal{A}\boldsymbol{x}_{k}^{m-1}) $} with PJ, PGS, and PF is $n\text{, } n^{2}$, and $n^{3}$, respectively.
Additionally, the J3 method \cite{r6} is the same as the J1 method \cite{r3} in program implementation.
Because the TAAR method with PF has \textcolor{black}{higher} computational cost each iteration than the TAAR method with the other two preconditioners, we compare the main computational cost each iteration in J1, GS1, J1\_SORlike, GS1\_SORlike, J2, GS2, GS3, FULLM, and TAAR with PF.

Let $ \boldsymbol{x}\text{, } \boldsymbol{y} \in \mathbb{R}^{n} $,  $ A\in \mathbb{R}^{n \times n} $, $ \mathcal{A} \in \mathbb{R}^{[m\text{,}n]} $.
Integers denote the frequency of each operation in the current method,
``$-$" means that the current algorithm does not have \textcolor{black}{that} operation.

Table \ref{t1} shows that the total computational cost of every method is $ O\left(mn^{m}\right) $, which indicates that it is feasible to compare the TAAR method with the other eight methods.

\begin{table}[H]
	\centering
	\caption{The main computational cost of nine methods per iteration.}
	\label{t1}
	\resizebox{.7\textwidth}{!}{
		\begin{tabular}{llllll}
			\hline
			&           & $ \mathcal{A} \boldsymbol{x} ^{m-1}   $ & 	$ \mathcal{A} \boldsymbol{x} ^{m-2} $ & $	 A\boldsymbol{x}=\boldsymbol{y}    $  & Total computational cost \\ \hline
			\multicolumn{2}{l}{J1}  &  $2$ &  $-$  & $-$   &   $O\left(mn^{m}\right)    $  \\
			\multicolumn{2}{l}{GS1}  & $2$  & $-$  &  $-$ &  $O\left(mn^{m}\right)   $  \\
			\multicolumn{2}{l}{J1\_SORlike}  &  $2$ &  $-$  & $-$   & $O\left(mn^{m}\right)   $  \\
			\multicolumn{2}{l}{GS1\_SORlike}  &$2$  & $-$  &  $-$  & $O\left(mn^{m}\right)   $   \\
			\multicolumn{2}{l}{J2}  & $1$  & $1$  &  $-$ &  $O\left(mn^{m}\right)   $  \\
			\multicolumn{2}{l}{GS2}  &$1$  & $1$  &  $-$ &  $O\left(mn^{m}\right)   $  \\
			\multicolumn{2}{l}{GS3}  &  $2$ &  $-$  & $-$ &  $O\left(mn^{m}\right)   $  \\
			\multicolumn{2}{l}{FULLM}  &   $2$ &  $-$  & $1$& $O\left(mn^{m}\right)   $   \\  \hline
			\multicolumn{1}{l|}{\multirow{2}{*}{TAAR with PF}}      & TR        &  $-$ &  $1$ &   $1$   & $O\left(mn^{m}\right)   $   \\ 
			\multicolumn{1}{l|}{}				& TAR         &  $-$ &  $1$ &   $1$  & $O\left(mn^{m}\right)   $   \\ \hline
			\multicolumn{2}{l}{Operation cost} & $ mn^{m}-n$   &  $\left(m-1\right)n^{m} - n^{2}   $   & $ n^{3}  $   & $-$  \\  \hline
	\end{tabular}}
\end{table}

\section{Numerical experiments} 
\label{sec:sec4}
We compare the TAAR method with three tensor splitting methods \cite{r3, r5, r6}.
\textcolor{black}{According to Section \ref{subsubsec:subsubsec3.2.1}, the TAAR method has three preconditioners: PJ, PGS, and PF.
	As mentioned in Section \ref{subsec:sec2.2}, the tensor splitting method 1 includes J1, GS1, J1\_SORlike, GS1\_SORlike methods. The tensor splitting method 2 includes J2 and GS2 methods. The tensor splitting method 3 includes J3, GS3 and FULLM methods.}
All experiments were performed in MATLAB R2021a, with the configuration: Inter(R) Core(TM) i7-10875H CPU at 2.30GHz CPU and 16.00GB RAM. 
We used the tensor toolbox 3.2.1 \cite{r27} to generate tensors and compute tensor products. We \textcolor{black}{used} Moore--Penrose pseudoinverse to compute $ \left(R_{k}^{\top}R_{k}\right)^{-1} $. 
The parameters $ \{ p\text{, }q \} $ \textcolor{black}{were} chosen as $ \{ p\text{, }q \}  =\{ 10\text{, }6 \} $ according to \cite{r16}.
We set the initial vector as $ \boldsymbol{x}_{0}=[0.1\text{, }0.1\text{,}\dots\text{, }0.1]^{\top}  $, 
the maximum number of iterations as 20,000
and the stopping criterion as
$$ \frac{\Vert \boldsymbol{b}-\mathcal{A} \boldsymbol{x}_{k}^{m-1} \Vert _{2} }{\Vert \boldsymbol{b}-\mathcal{A} \boldsymbol{x}_{0}^{m-1} \Vert _{2} }  \leqslant 10^{-8}\text{.} $$

\begin{emp}
	\label{emp1}
\end{emp}
According to \cite{r3}, we construct a nonsingular $ \mathcal{M} $ tensor $ \mathcal{A} = s \mathcal{I} - \mathcal{B} \in \mathbb{R}^{\left[m\text{,}n\right]} $ with $ s $ satisfying
$$ s=\left(1+\varepsilon\right) \max _{i=1\text{,}\dots\text{, }n} \left(\mathcal{B} \boldsymbol{e}^{n-1}\right)_{i}\text{,} \quad  \varepsilon > 0 \text{,}$$
where $ \boldsymbol{e} = \left(1\text{, }1\text{,}\dots\text{, }1\right)^{\top} $ and $ \varepsilon = 0.01 $.
We choose ``rand'' as rand(`state', $0$).
The entries of $ \mathcal{B} $ are generated from the standard uniform distribution on $ \left(0\text{, }1\right) $ by ``tenrand" 
and positive $ \boldsymbol{b} $ is generated from the same distribution by ``rand".
Obviously, $ \mathcal{A} $ is a $ \mathcal{Z} $-tensor, which satisfies $  \mathcal{A} \boldsymbol{e}^{n-1}>0 $. According to Proposition \ref{lma1}, $ \mathcal{A} $ is a nonsingular $ \mathcal{M} $-tensor.

We \textcolor{black}{perform} seven groups of $ \left(m\text{, }n\right) $ and randomly generated the corresponding \textcolor{black}{datasets} $ \left(\mathcal{A}\text{, }\boldsymbol{b}\right) $. The cases of $ \left(m\text{, }n\right)$ are $ \left(m\text{, }n\right) = \{\left(3\text{, }200\right)\text{,} \left(3\text{, }400\right)\text{,} \left(3\text{, }600\right)\text{,} \left(4\text{, }50\right)\text{,} \left(4\text{, }100\right)\text{,} \left(5\text{, }20\right)\text{,} \left(5\text{, }40\right)\}$.

We \textcolor{black}{test} the effect of choosing different preconditioners on the convergence of the TAAR method.
As mentioned in Section $\ref{subsubsec:subsubsec3.2.1}$,  
the TAAR method has three preconditioners, respectively, labeled as PJ, PGS, and PF. 
We \textcolor{black}{use} the number of iterations, the normalized residual, and the CPU time in seconds to measure the convergence performance, respectively, denoted by Iter, Res, and CPU[s]. 

The stopping criterion for the \textcolor{black}{experiment} in Table \ref{t2} is the same. We set the stopping criterion as the relative residual less than or equal to $ 10^{-8} $. Observing one case $ \left(m\text{, }n\right) =\left(4\text{, }50\right)$ in Table \ref{t2}, the number of iterations of TAAR with PJ is $ 19 $.  The TAAR method with PJ does not satisfy the stopping criterion at the $18$th iteration, but satisfies the stopping criterion at the $19$th iteration and its normalized residual is $ 6.291 \times 10^{-12} $. It indicates that the relative residual of TAAR with PJ decreases rapidly from the $18$th iteration to the $19$th iteration.

Observing one case $ \left(m\text{, }n\right) =\left(3\text{, }200\right)$ in Table \ref{t2}, the number of iterations, normalized residual, CPU time of the TAAR method with PF are similar to those of the TAAR method with the other two preconditioners, and this conclusion is also suitable for other cases. We \textcolor{black}{conclude} that the performance of the TAAR method with different preconditioners is almost the same.

Because the TAAR method with PF has \textcolor{black}{higher} computational cost than the TAAR method with \textcolor{black}{the} other two preconditoners, the TAAR method with PF may require more CPU time.
Based on the above conclusion, we chose the TAAR method with PF as an example to compare with other existing tensor splitting methods in the following experiments.

\begin{table}[H]
	\caption{The comparison of the TAAR method with three preconditioners.}
	\label{t2}
	\setlength\tabcolsep{5pt}
	\centering
	\resizebox{.9\textwidth}{!}{
		\begin{tabular}{llllllllll}
			\hline
			\multirow{2}{*}{(m\text{, }n)} & \multicolumn{3}{l}{Iter} & \multicolumn{3}{l}{Res} & \multicolumn{3}{l}{CPU[s]} \\ 
			\cmidrule(lr){2-4} \cmidrule(lr){5-7} \cmidrule(lr){8-10}
			&PJ       & PGS      & PF     &PJ       & PGS      & PF     &PJ       & PGS      & PF              \\ \hline
			$\left(3\text{, }200\right)$ & $16$ & $15$ & $18$ & $8.950 \times 10^{-9} $  & $5.938 \times 10^{-9} $ & $9.551 \times 10^{-9} $  & $0.039$   &$ 0.035$    &  $0.046$     \\
			$\left(3\text{, }400\right)$ & $13$  & $13$ &$ 13$ & $3.940 \times 10^{-10} $  & $1.738 \times 10^{-9} $  & $3.889 \times 10^{-10} $  & $0.246$  & $0.232$ &$  0.252$    \\
			$\left(3\text{, }600\right)$ &$13 $& $13$&$ 13$     & $6.186 \times 10^{-9} $  & $3.727 \times 10^{-9} $ & $5.326 \times 10^{-9} $& $0.833$   & $0.766$    &$ 0.797$  \\
			$\left(4\text{, }50\right)$  &$19$  & $17$  & $16$    & $6.291 \times 10^{-12} $ & $4.620 \times 10^{-9} $  & $4.646 \times 10^{-9} $  & $0.038$   &$0.027$     & $0.026$   \\
			$\left(4\text{, }100\right)$  & $14$    & $15$   & $14$  & $5.192\times 10^{-9} $ & $2.430 \times 10^{-9} $ & $8.820 \times 10^{-9} $   & $0.478 $  & $ 0.580$    & $0.500$  \\ 
			$\left(5\text{, }20\right)$  & $19$   & $19$   & $18$    & $7.796 \times 10^{-11} $  & $1.108 \times 10^{-11} $  & $6.405 \times 10^{-9} $   & $0.021$  & $0.023$  & $0.026$     \\
			$\left(5\text{, }40\right)$  & $19$     & $19$      & $16$    & $1.803\times 10^{-12} $   & $1.525 \times 10^{-12} $   & $9.807 \times 10^{-9} $      & $0.881$   & $0.853$   & $0.719$   \\ \hline
	\end{tabular}}
	\end{table}
	
	\begin{emp}
\label{emp2}
\end{emp}

In this experiment, we \textcolor{black}{compare} the computational \textcolor{black}{cost} between the TAAR method with PF and three tensor splitting methods in Section \ref{subsec:subsec2.2}.
\textcolor{black}{We select three cases from Experiment \ref{emp1}. The cases are $ \left(m\text{, }n\right) = \{\left(3\text{, }200\right)\text{,} $ $\left(3\text{, }400\right)\text{,} $ $\left(4\text{, }100\right) \}$.
Since the CPU time of the GS1 and GS1\_SORlike methods in real computation was more than half an hour, we only compared the TAAR method with the J1, J1\_SORlike, J2, GS2, GS3, and FULLM methods.}

\begin{figure}[H]
\centering
\subfigure{
	\begin{minipage}[h]{0.315\linewidth}
		\centering
		\includegraphics[width=1\linewidth]{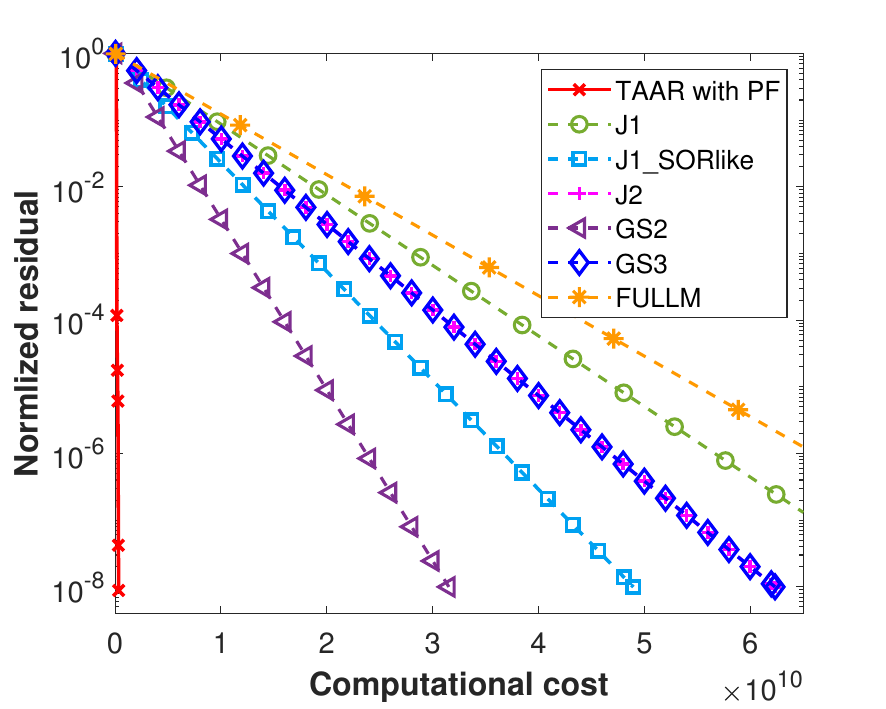}
		\caption*{(a) $\left(m\text{, }n\right)=\left(3\text{, }200\right)$}
	\end{minipage}
}
\subfigure{
	\begin{minipage}[h]{0.315\linewidth}
		\centering
		\includegraphics[width=1\linewidth]{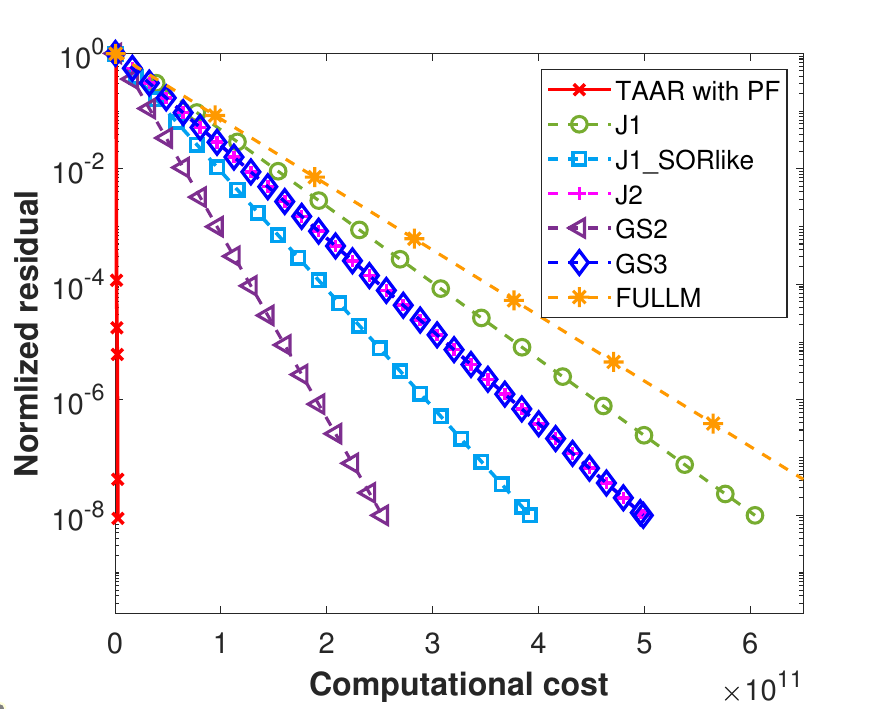}
		\caption*{(b) $\left(m\text{, }n\right)=\left(3\text{, }400\right)$}
	\end{minipage}
}
\subfigure{
	\begin{minipage}[h]{0.315\linewidth}
		\centering
		\includegraphics[width=1\linewidth]{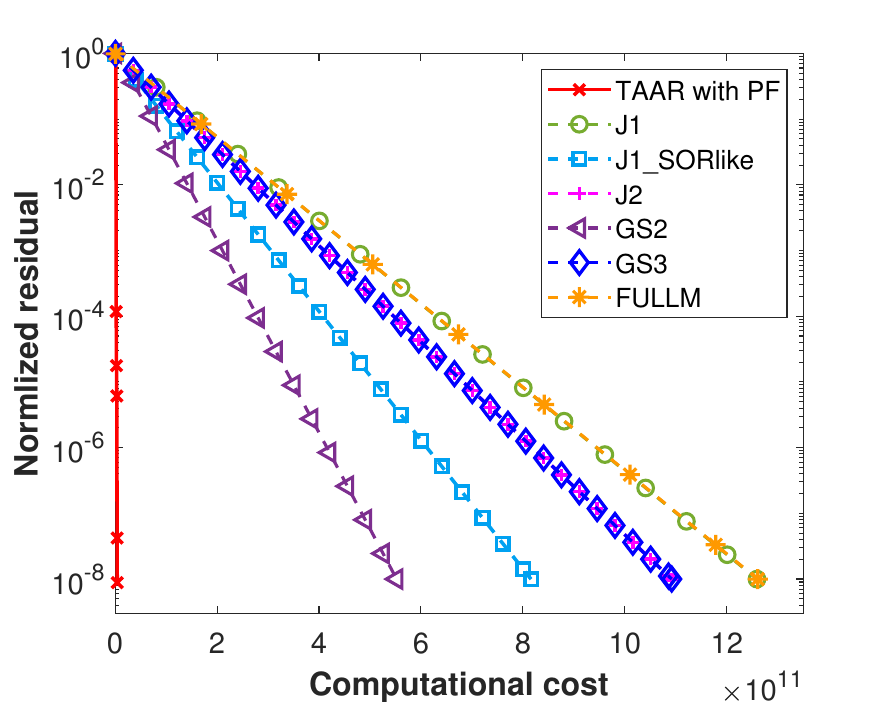}
		\caption*{(c) $\left(m\text{, }n\right)=\left(4\text{, }100\right)$}
	\end{minipage}
}
\setlength{\abovecaptionskip}{1pt}
\caption{The comparison of computational \textcolor{black}{cost} between the TAAR method with PF and other splitting methods. \protect\\ $x$-axis: the sum of flops up to the current iteration step.}
\label{fig1}
\end{figure}

When achieving the same relative residual accuracy,
Fig. \ref{fig1} shows that the computational cost of the TAAR method is lower than those of the tensor splitting methods; the GS2 method had the second lowest computational cost, followed by the J2 method, and the computational \textcolor{black}{cost} of the J1, J2, and GS3 methods was similar.
From Fig. \ref{fig1},
we can conclude that the TAAR method is convergent.
In addition, the convergence of the TAAR method seems to be linear, which indicates that the technique we used to compute the relaxation parameters $\omega_{k}\text{ and }\beta_{k}$ is valid.

\begin{emp}
\label{emp3}
\end{emp}

Based on the conclusion of Experiment \ref{emp2}, we further \textcolor{black}{test} the effectiveness of the TAAR method with PF in seven cases from Experiment \ref{emp1}.
We \textcolor{black}{compare} the TAAR method with the J1, J1\_SORlike, J2, GS2, GS3, and FULLM methods.


From Table \ref{t3}, the TAAR method \textcolor{black}{is} convergent and \textcolor{black}{can} achieve the same accuracy within fewer iterations and shorter time than other splitting methods. 
Table \ref{t3} shows that the TAAR method can accelerate the convergence by at least one order of magnitude compared with the existing tensor splitting methods. It indicates that the effectiveness of applying Anderson acceleration to the tensor splitting method \cite{r6} for solving Eq. \eqref{Eq1.1}.

\begin{table}[H]
\caption{The comparison of the proposed method with three tensor splitting methods.}
\label{t3}
\setlength\tabcolsep{5pt}
\centering
\resizebox{.9\textwidth}{!}{
	\begin{tabular}{lllllllll}
		\hline
		& & $\left(3\text{, }200\right)$ & $\left(3\text{, }400\right)$ & $\left(3\text{, }600\right)$ & $\left(4\text{, }50\right)$ & $\left(4\text{, }100\right)$ & $\left(5\text{, }20\right)$ & $\left(5\text{, }40\right)$  \\  \hline
		\multirow{3}{*}{TAAR with PF (proposed)} &  Iter &  $18$ & $ 13 $ & $ 13 $ & $ 16 $ & $ 14 $ &  $ 18 $&  $ 16 $ \\
		&  Res &   $9.551 \times 10^{-9}$ & $ 3.889 \times 10^{-10} $  &  $ 5.326 \times 10^{-9} $ & $ 4.646 \times 10^{-9} $  & $ 8.820\times 10^{-9}  $  & $6.405 \times 10^{-9}   $  & $ 9.807 \times 10^{-9}  $  \\
		& CPU[s]  & $ 0.050 $  &  $ 0.291 $ &$ 0.869 $   & $ 0.033  $  & $ 0.522 $  & $ 0.034 $  & $ 0.808 $  \\  \hline
		\multirow{3}{*}{J1} &  Iter & $ 1054 $  & $ 1335 $  &  $ 1413 $ &  $ 1234 $ & $ 1574 $  & $ 1415 $  & $ 1697 $  \\
		& Res  &  $ 9.886 \times 10^{-9} $ & $ 9.962 \times 10^{-9} $  & $ 9.952\times 10^{ -9}  $ & $ 9.998 \times 10^{-9} $ &$ 9.961\times 10^{-9}  $ & $ 9.912 \times 10^{-9}  $& $ 9.985\times 10^{-9}  $ \\
		& CPU[s]   &  $ 4.232 $  & $ 43.419 $& $ 152.338 $& $ 3.944 $ & $ 95.686 $ &$ 3.282 $ & $ 143.073 $ \\  \hline
		\multirow{3}{*}{J1\_SORlike} & Iter  &  $682$ & $ 865 $ & $ 915 $ & $ 799 $ &$ 1020 $ & $ 917 $& $ 1100 $ \\
		& Res  & $ 9.849 \times 10^{-9} $  &$ 9.859 \times 10^{-9} $& $ 9.994\times 10^{-9}  $ & $ 9.978\times 10^{-9}  $ & $ 9.936 \times 10^{-9}  $ &$ 9.823 \times 10^{-9}  $&$ 9.954 \times 10^{-9}  $  \\
		& CPU[s]  & $2.658$  & $ 28.055 $ &$ 99.942 $  & $ 2.468 $ & $ 62.008 $ &$ 2.152 $ &  $ 92.679 $  \\  \hline
		\multirow{3}{*}{J2} &  Iter &  $1050$ & $ 1334 $  & $ 1413 $ & $ 1211 $ & $ 1560 $ &$ 1345 $ &$ 1657 $ \\
		& Res  & $ 9.930 \times 10^{-9} $  &$ 9.890 \times 10^{ -9} $  & $ 9.904 \times 10^{-9}  $ & $ 9.954\times 10^{-9}  $ &$ 9.981\times 10^{-9} $  & $ 9.904\times 10^{-9} $&$ 9.899\times 10^{-9} $ \\
		& CPU[s]   &  $ 4.059 $ & $ 43.675 $& $ 153.184 $ & $ 3.709 $ & $ 95.630 $& $ 3.154 $ & $ 137.535 $\\  \hline
		\multirow{3}{*}{GS2} & Iter  & $ 532 $  &$ 675 $  & $ 715 $ & $ 613 $ &$ 789 $ & $ 680 $&$  838$ \\
		& Res  & $ 9.842 \times 10^{-9} $  &$ 9.883\times 10^{ -9}  $ & $ 9.874\times 10^{-9 }  $ & $ 9.870\times 10^{-9}  $ &$ 9.971\times 10^{-9} $ &  $ 9.931\times 10^{-9} $& $ 9.816\times 10^{-9} $ \\
		& CPU[s]   & $ 2.080 $  & $ 22.224 $&$ 76.914 $  &$ 1.872 $ &  $ 51.507 $& $ 1.576 $ &$ 73.845 $  \\  \hline
		\multirow{3}{*}{GS3} & Iter  & $ 1051 $  & $ 1334 $ & $ 1412 $ &$ 1234 $  & $ 1574 $ &$ 1415 $ & $ 1697 $\\
		&  Res & $ 9.955 \times 10^{-9} $  &$ 9.872\times 10^{-9}  $  &$ 9.931\times 10^{-9}  $ &$  9.963\times 10^{-9}  $ & $ 9.953 \times 10^{-9} $ &$ 9.902\times 10^{-9} $  & $ 9.983\times 10^{-9} $ \\
		& CPU[s]  & $ 4.157 $  & $ 43.522 $ & $ 153.564 $ &$ 3.835 $  &$ 104.124 $  & $ 3.328 $ &$ 143.188 $ \\  \hline
		\multirow{3}{*}{FULLM} & Iter  & $ 1049 $  &  $ 1332 $& $ 1411 $ & $ 1234 $ & $ 1574 $&  $ 1415 $& $ 1697 $ \\
		& Res  & $ 9.851 \times 10^{-9} $  & $ 9.918\times 10^{-9}  $ &$ 9.909\times 10^{-9}  $  &  $ 9.928\times 10^{-9}  $ & $ 9.943\times 10^{-9} $& $ 9.891\times 10^{-9} $ &  $ 9.982\times 10^{-9} $\\
		&  CPU[s]  &  $ 4.369 $ & $ 45.217 $ &$ 155.323 $ & $ 3.905 $&$ 113.944 $ & $ 3.223 $&$ 143.077 $  \\  \hline
	\end{tabular}
}
\end{table}

\textcolor{black}{\begin{emp}
	\label{emp4.0}
\end{emp}
This experiment is from \cite{r9}. Let $ s=n^2 $. We construct a $3$rd-order symmetric $\mathcal{M}$-tensor $ \mathcal{A}=s \mathcal{I}-\mathcal{B} $ with
$$ b_{ijk}= \left| \sin (i+j+k) \right| \text{. }$$
The right-hand side $\boldsymbol{b}= (1\text{, }1\text{, }\dots\text{, }1)^{\top} $.
We compare the proposed method with the Newton method in \cite{r3}.
Table \ref{t4.0} shows the proposed method is more efficient than the Newton method if $ n$ is large. }

\begin{table}[H]
\captionsetup{labelfont={color=black}}
\arrayrulecolor{black}
\caption{\color{black}{The comparison of TAAR with PF and Newton method.}}
\label{t4.0}
\setlength\tabcolsep{6pt}
\centering
\resizebox{.7\textwidth}{!}{
	\begin{tabular}{lllllll}
		\hline
		\multirow{2}{*}{\textcolor{black}{(m, n)}} & \multicolumn{3}{l}{\textcolor{black}{TAAR with PF (proposed)}} & \multicolumn{3}{l}{\textcolor{black}{Newton method}}   \\ 
		\cmidrule(lr){2-4} \cmidrule(lr){5-7}
		&\textcolor{black}{Iter}      & \textcolor{black}{Res}      & \textcolor{black}{CPU[s]}     &\textcolor{black}{Iter}      & \textcolor{black}{Res}      & \textcolor{black}{CPU[s]}                  \\ \hline
		\textcolor{black}{$\left(3\text{, }50\right)$} & \textcolor{black}{$7$} & \textcolor{black}{$2.450\times 10^{-12}$} & \textcolor{black}{$0.001$} & \textcolor{black}{$5 $}  & \textcolor{black}{$1.288 \times 10^{-10} $} & \textcolor{black}{$0.002 $ }     \\
		\textcolor{black}{$\left(3\text{, }100\right)$} & \textcolor{black}{$6$} & \textcolor{black}{$7.878\times 10^{-9}$} & \textcolor{black}{$0.002$} & \textcolor{black}{$6 $}  & \textcolor{black}{$5.516 \times 10^{-11} $} & \textcolor{black}{$0.005 $   }      \\
		\textcolor{black}{$\left(3\text{, }200\right)$ }& \textcolor{black}{$7$} & \textcolor{black}{$3.010\times 10^{-12}$ }& \textcolor{black}{$0.019$} & \textcolor{black}{$7 $ } & \textcolor{black}{$1.587 \times 10^{-11} $} & \textcolor{black}{$0.034 $ }     \\
		\textcolor{black}{$\left(3\text{, }300\right)$} & \textcolor{black}{$6$} & \textcolor{black}{$9.892\times 10^{-9}$} & \textcolor{black}{$0.061$} & \textcolor{black}{$7 $ } & \textcolor{black}{$8.604 \times 10^{-9} $ }& \textcolor{black}{$0.115 $ }     \\
		\textcolor{black}{$\left(3\text{, }400\right)$} & \textcolor{black}{$4$} & \textcolor{black}{$5.693\times 10^{-9}$} & \textcolor{black}{$0.084$} & \textcolor{black}{$8 $}  & \textcolor{black}{$4.086 \times 10^{-12} $} & \textcolor{black}{$0.275 $}  			  \\ \hline
\end{tabular}}
\end{table}

\begin{emp}
\label{emp4}
\end{emp}

\begin{figure}[H]
\centering
\includegraphics[scale=0.4]{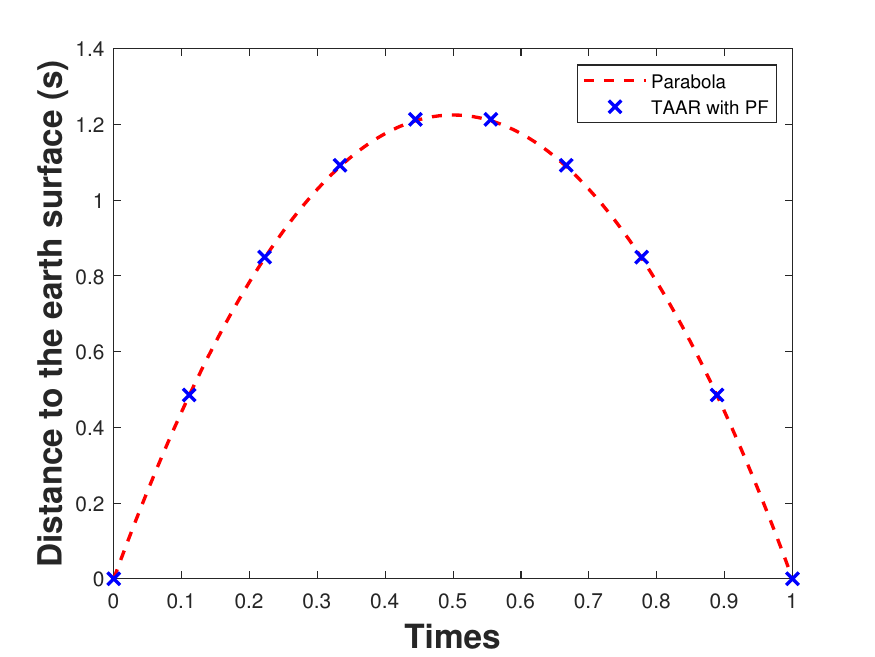}
\setlength{\abovecaptionskip}{8pt}
\caption{\textcolor{black}{The results for Experiment \ref{emp4}}.}
\label{fig4}
\end{figure}

This experiment comes from \cite{r3}. Consider the ordinary differential equation
$$ \frac{\mathrm{d}^{2} x(t)}{\mathrm{d} t^{2}} = -\frac{GM}{x(t)^{2}}\text{,} \quad t\in(0\text{, }1) $$
with Dirichlet's boundary conditions
$$ x(0) = c_{0} \text{,} \quad  x(1) = c_{1} \text{. } $$
The above equation can describe a particle's movement under the gravitation
\textcolor{black}{$$ m\frac{\mathrm{d}^{2} x(t)}{\mathrm{d} t^{2}} = -\frac{GMm}{x(t)^{2}}\text{, }$$
where $ G  \approx 6.67\times 10^{-11} \text{Nm}^{2}/\text{kg}^{2}$ is the gravitational constant and $ M \approx 5.98 \times 10^{24}\text{kg} $ is the mass of the earth.
Assuming that the distance between the earth’s surface and the earth’s center is $ 6.37\times 10^{6} $ meters, we consider the trajectory of a particle after it is thrown upward near the earth’s surface. The trajectory can be approximated by a parabola
\begin{equation} \nonumber
\left\{
\begin{aligned}
	& x_{t}= -\frac{1}{2}gt^{2}+\alpha+\beta \text{, } \\
	& x(0)=c_{0}\text{,} \quad  x(1)=c_{1} \text{, }
\end{aligned}
\right.
\end{equation}
where $g\approx 9.8 \text{m}/ \text{s}^{2}$ and $c_{0}=c_{1}=6.37\times 10^{6}$.
We plot the trajectory in Fig. \ref{fig4}.}

After discretizing the ordinary differential equation, we can rewrite the discretization into \textcolor{black}{a} multilinear system,
$$ \mathcal{A} x^{3} =b \text{, }$$
where $\mathcal{M}$-tensor $ \mathcal{A} \in \mathbb{R}^{[4,n]}$ satisfies
\begin{equation} \nonumber
\left\{
\begin{aligned}
& a_{1111}=a_{nnnn}=1\text{, } \\
& a_{iiii}=2 \text{,}\quad  i =2\text{, } 3 \text{,} \dots \text{, } n-1\text{, } \\
&    a_{i(i-1)ii} = a_{ii(i-1)i} = a_{iii(i-1)} = -\frac{1}{3} \text{,} \quad  i =2\text{, } 3 \text{,} \dots \text{, } n-1\text{, } \\
& a_{i(i+1)ii} = a_{ii(i+1)i} = a_{iii(i+1)} = -\frac{1}{3} \text{,}\quad  i =2\text{, } 3 \text{,} \dots \text{, } n-1\text{, } 
\end{aligned}
\right.
\end{equation}
and $b$ satisfies
\begin{equation} \nonumber
\left\{
\begin{aligned}
& b_{1}=c^{3}_{0}\text{, } \\
& b_{i}=\frac{GM}{(n-1)^{2}} \text{,} \quad  i =2\text{, } 3 \text{,} \dots \text{, } n-1\text{, } \\
& b_{n}=c^{3}_{1} \text{. }
\end{aligned}
\right.
\end{equation}
\textcolor{black}{We solve the multilinear system by using the TAAR method with PF, and set $ n=20 $. Figure \ref{fig4} shows that the solution we obtained satisfies the real world.}

\section{Concluding remarks} 
\label{sec:sec5}
Inspired by the AAR method for solving Eq. \eqref{Eq3.1}, we proposed a TAAR method for solving Eq. \eqref{Eq1.1}.
We first presented a TR method based on tensor regular splittings, then applied Anderson acceleration to the TR method and derived a TAR method, and finally, we proposed \textcolor{black}{a} TAAR method by periodically employing the TAR method within the TR method.
Numerical experiments showed that the TAAR method could accelerate convergence by at least one order of magnitude compared with other existing tensor splitting methods. In addition, the TAAR method could achieve the same accuracy within fewer iteration numbers and \textcolor{black}{a} shorter time than other tensor splitting methods.
\textcolor{black}{In addition, there are many structure tensors such as $\mathcal{M}$-tensor, $\mathcal{H}$-tensor, $\mathcal{L}$-tensor, $\mathcal{Z}$-tensor. And there are many algorithms for solving the multilinear systems with these special structure tensors. Studies on these special structures of tensors will be interesting future work.}

\section*{Appendix}
\begin{figure}[H]
	\centering
	\includegraphics[scale=0.4]{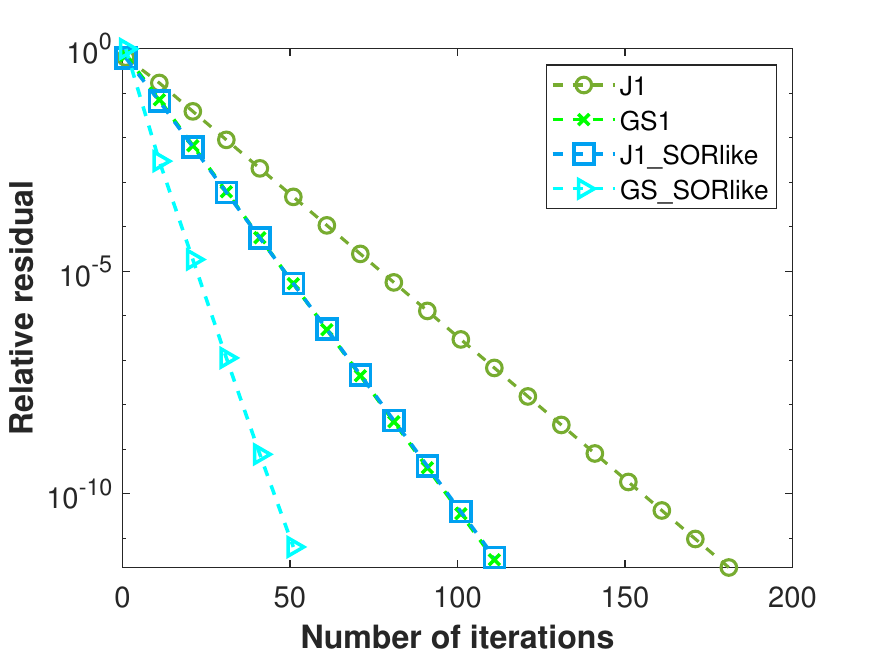}
	\setlength{\abovecaptionskip}{8pt}
	\captionsetup{labelfont={color=black}}
	\caption{\color{black}{The results for Experiment \ref{emp5}}.}
	\label{fig5}
\end{figure}

\textcolor{black}{We compare the implementations of tensor splitting methods 1, 2, and 3 with the prior works \cite{r3, r5, r6}. The CPU time is different from that in prior works because of different environment configurations.}
\textcolor{black}{	
	\begin{emp}
		\label{emp5}
	\end{emp}
	We compare the implementations of the tensor splitting method 1 with the prior work \cite{r3}.
	We generate a $3$rd-order $10$-dimensional nonsingular $ \mathcal{M} $-tensor $ \mathcal{A} $ by using the Experiment \ref{emp1} and set $  \varepsilon = 0.01$.
	The acceleration parameter is $\omega  = 0.35 \cdot \min _{i = 1 \text{, }2 \text{, } \dots \text{, } n} a_{ii \dots i}  $.
	The way to choose a right-hand side $ \boldsymbol{b} $ and an initial value are not given in \cite{r3}, so it's hard to get the same number of iterations as the Example 4.1 in \cite{r3}. We reset the generator to the $1$st state by using ``rand(`state', 1)" in this experiment. The stopping criterion is
	$$ \frac{\Vert \boldsymbol{b}-\mathcal{A} \boldsymbol{x}_{k}^{m-1} \Vert _{2} }{\Vert \boldsymbol{b}-\mathcal{A} \boldsymbol{x}_{0}^{m-1} \Vert _{2} }  \leqslant 10^{-12}\text{.} $$
	Figure \ref{fig5} shows that the convergence of these four methods in tensor splitting method 1 is close to the convergence of Fig. 3 in \cite{r3}.}

\textcolor{black}{\begin{emp}
		\label{emp6}
	\end{emp}
	We compare the implementations for the tensor splitting method 2 with the prior work \cite{r5}.
	This experiment is the same as Problem 3 from \cite{r5}.
	Set $\mathbf{b}=(1\text{, }2\text{, }3)^{\top}$, $ \mathbf{x}_{0} = (1 \text{, } 1 \text{,} 1)^{\top} $, and a symmetric tensor $ \mathcal{A} \in \mathbb{R}^{[4,3]}$, elementwise, $a_{1111}=20.4982$, $a_{1112} = −0.0582$, $a_{1113} = −1.1719$, $a_{1122} =0.2236$, $a_{1123} = −0.0171$, $a_{1133} = 0.4597$, $a_{1223} = 0.1852$, $a_{1222} = 0.4880$, $a_{1233} = −0.4087$, $a_{1333} = 0.7639$, $a_{2222} = 10$, $a_{2223} = −0.6162$, $a_{2233} = 0.1519$, $a_{3333} = 2.6311$. The initial vector is $ \mathbf{x}_{0}= (1 \text{, }1 \text{, } 1)^{\top}$. We denote J2, GS2 in \cite{r5} as J2 (prior), GS2 (prior). Table \ref{te5} shows that the number of iterations of J2 and GS2 in the tensor splitting method 2 are the same as those in \cite{r5}.}
\begin{table}[H]
	\setlength{\belowcaptionskip}{8pt}
	\captionsetup{labelfont={color=black}}
	\arrayrulecolor{black}
	\caption{\color{black}{Numerical results for Experiment 5.2.}\label{te5}}
	\centering
	\resizebox{.35\textwidth}{!}{
		\begin{tabular}{lll}
			\hline
			& \textcolor{black}{Iter}  & \textcolor{black}{CPU[s]}  \\ \hline
			\textcolor{black}{J2 (our implementation)}    & \textcolor{black}{$14$} & \textcolor{black}{0.0013}  \\ \hline
			\textcolor{black}{J2 (prior work) }   & \textcolor{black}{$14$} & \textcolor{black}{0.0010}  \\ \hline
			\textcolor{black}{GS2 (our implementation)} & \textcolor{black}{$10$}  & \textcolor{black}{0.0007}  \\ \hline
			\textcolor{black}{GS2 (prior work)} & \textcolor{black}{$10$}  & \textcolor{black}{0.0007}  \\ \hline
		\end{tabular}
	}
\end{table}

\textcolor{black}{\begin{emp}
		\label{emp7}
	\end{emp}
	We use Experiment \ref{emp1} to compare the implementations of the tensor splitting method 3 with the prior work \cite{r6}.
	This Experiment setting is the same as Example $6.1$ from \cite{r6}.
	We consider a $3 $rd-order $5$-dimensional nonsingular $  \mathcal{M}$-tensor, and set $ \epsilon=1 $ and $ \boldsymbol{b} = \boldsymbol{x}_{0} = (1 \text{, } 1 \text{,} \dots \text{, } 1)^{\top} $. The stopping criterion is
	$$ \Vert \boldsymbol{b}-\mathcal{A} \boldsymbol{x}_{k}^{m-1} \Vert _{2}  \leqslant 10^{-11}\text{.} $$
	Table \ref{te6} shows that the implementations of the tensor splitting method 3 are as effective as the prior works \cite{r6}.}

\begin{table}[H]
	\setlength{\belowcaptionskip}{8pt}
	\arrayrulecolor{black}
	\captionsetup{labelfont={color=black}}
	\caption{\color{black}{Numerical results for Experiment 5.3.}\label{te6}}
	\centering
	\resizebox{.35\textwidth}{!}{
		\begin{tabular}{lll}
			\hline
			& \textcolor{black}{Iter}  & \textcolor{black}{CPU[s]}  \\ \hline
			\textcolor{black}{GS2 (our implementation)}    & \textcolor{black}{$33$} & \textcolor{black}{0.0023}  \\ \hline
			\textcolor{black}{GS2 (prior work)}   & \textcolor{black}{$34$}  & \textcolor{black}{0.0027}  \\ \hline
			\textcolor{black}{FULLM2 (our implementation)} & \textcolor{black}{$31$}  & \textcolor{black}{0.0035}  \\ \hline
			\textcolor{black}{FULLM2 (prior work)} & \textcolor{black}{$31$}  & \textcolor{black}{0.0065}  \\ \hline
		\end{tabular}
	}
\end{table}

\section*{Declarations}

\subsection*{Competing interests}

The authors declare no competing interests.

\subsection*{Generative AI in scientific writing}

The authors declare that no AI was used in the writing process.

\subsection*{Acknowledgments}
The authors appreciate the anonymous reviewers for their fruitful comments that enhanced the quality of the manuscript.
This work was supported by the China Scholarship Council and supported by JSPS KAKENHI Grant Number: JP20H00581.

\bibliographystyle{unsrt}
\bibliography{R1}

\end{document}